\documentclass[11pt]{amsart}
\usepackage{amssymb, adjustbox, enumerate, amsbsy, stmaryrd}
\usepackage{geometry}
\usepackage{todonotes}

\usepackage{amsfonts, amssymb, amscd}
\numberwithin{equation}{section}

\usepackage[symbol]{footmisc}

\usepackage{bm}
\usepackage{verbatim}
\usepackage{mathrsfs}
\usepackage{graphicx}
\usepackage{tikz-cd}
\usepackage{subcaption}
\usepackage{listings}
\usepackage{subfiles}
\usepackage[toc,page]{appendix}
\usepackage{mathtools}
\usepackage{comment}
\usepackage{enumerate}
\usepackage{enumitem}
\usepackage[all]{xy}

\usepackage{graphicx}
\graphicspath{{images/}}

\usepackage{appendix}
\usepackage{hyperref}
\hypersetup{
    colorlinks=true,
    citecolor=red,
    linkcolor=blue,
    filecolor=magenta,      
    urlcolor=red,
}
\newcommand{\bb}{\bm{b}}
\newcommand{\Mm}{{\bf{M}}}
\newcommand{\Bb}{{\bf{B}}}
\newcommand{\PP}{{\bf{P}}}
\newcommand{\NN}{{\bf{N}}}

\newcommand{\Spec}{\mathrm{Spec}}
\newcommand{\Src}{\mathrm{Src}}

\newcommand{\id}{\mathrm{id}}

\newcommand{\Cc}{\mathbb{C}}

\newcommand{\Pp}{\mathbb{P}}
\newcommand{\Qq}{\mathbb{Q}}

\newcommand{\Rr}{\mathbb{R}}

\newcommand{\Exc}{\operatorname{Exc}}

\newcommand{\Fix}{\operatorname{Fix}}
\newcommand{\Mov}{\operatorname{Mov}}
\newcommand{\Bir}{\operatorname{Bir}}
\newcommand{\Aut}{\operatorname{Aut}}

\newcommand{\rk}{\operatorname{rank}}
\newcommand{\red}{\operatorname{red}}

\newcommand{\lcg}{\operatorname{lcg}}
\newcommand{\Nklt}{\operatorname{Nklt}}
\newcommand{\mld}{{\rm{mld}}}

\newcommand{\Supp}{\operatorname{Supp}}
\newcommand{\Ngklt}{\operatorname{Ngklt}}
\newcommand{\Nlc}{\operatorname{Nlc}}

\newcommand{\mult}{\operatorname{mult}}

\newcommand{\lf}{\lfloor}
\newcommand{\rf}{\rfloor}

\newcommand{\Aa}{{\bf{A}}}

\newcommand{\Oo}{\mathcal{O}}

\newcommand{\Pic}{\mathrm{Pic}}

\newtheorem{thm}{Theorem}[section]

\newtheorem{cor}[thm]{Corollary}
\newtheorem{lem}[thm]{Lemma}

\theoremstyle{definition}
\newtheorem{defn}[thm]{Definition}
\newtheorem{ques}[thm]{Question}
\theoremstyle{definition}
\newtheorem{rem}[thm]{Remark}

\newtheorem{defthm}[thm]{Definition-Theorem}

\newtheorem{ex}[thm]{Example}

\newtheorem{cons}[thm]{Construction}

\theoremstyle{definition}

\begin{document}

\title{Semi-ampleness of NQC generalized log canonical pairs}
\author{Jihao Liu and Lingyao Xie}

\address{Department of Mathematics, Northwestern University, 2033 Sheridan Rd, Evanston, IL 60208, USA}
\email{jliu@northwestern.edu}

\address{Department of Mathematics, The University of Utah, Salt Lake City, UT 84112, USA}
\email{lingyao@math.utah.edu}

\date{\today}

\begin{abstract}
We establish a Koll\'ar-type gluing theory for NQC generalized log canonical pairs and use it to prove semi-ampleness results of NQC generalized pairs. As consequences, we prove the existence of flips for any NQC generalized log canonical pair, and show that NQC generalized log canonical singularities are Du Bois.
\end{abstract}

\maketitle
\tableofcontents

\section{Introduction}

We work over the field of complex numbers $\mathbb C$. We remark that we expect our results to hold over any algebraically closed field of characteristic zero. However, since many important references we cite only work over $\mathbb C$ (e.g. \cite{BZ16,HL21a,HL22}), this paper will also only work over $\mathbb C$ for consistency.

The theory of \emph{generalized pairs} (\emph{g-pairs} for short) is a central topic in modern day birational geometry. Introduced by Birkar and Zhang in \cite{BZ16} in the study of effective Iitaka fibrations, this theory is known to be useful in many aspects of birational geometry, such as the proof of the Borisov-Alexeev-Borisov conjecture \cite{Bir19,Bir21a}, the theory of complements \cite{Bir19,Sho20}, the connectedness principles \cite{Bir20,FS20}, non-vanishing theorems \cite{LMPTX22}, etc. We refer the reader to \cite{Bir21b} for a survey on the theory of g-pairs.

An important part of the study of g-pairs is their minimal model program. The foundations of the minimal model program for gklt g-pairs and $\mathbb Q$-factorial gdlt g-pairs were established in \cite{BZ16,HL22}. Recently, there is some progress towards the minimal model program theory for glc g-pairs. In particular, in \cite{HL21a}, the authors proved the cone theorem, contraction theorem, and the existence of flips for NQC $\Qq$-factorial glc g-pairs. For other related works, we refer the reader to \cite{Has20b,LT21,Has22,LX22}. These results almost complete the foundation of the minimal model program for $\mathbb Q$-factorial NQC glc g-pairs, or for NQC g-pairs admitting an lc structure on the ambient variety.

In this paper, we focus on the last part of the minimal model program for NQC generalized pairs: the class of possibly non-$\mathbb Q$-factorial NQC g-pairs. The main theorem of this paper is the following:
\begin{thm}\label{thm: gmm exists for g-crepant log structure}
Let $(X,B,\Mm)/Z$ be an NQC glc g-pair and $A\geq 0$ an $\Rr$-divisor on $X$, such that $(X,B+A,\Mm)$ is glc and $K_X+B+A+\Mm_X\sim_{\Rr,Z}0$. Then:
\begin{enumerate}
    \item $(X,B,\Mm)/Z$ has a Mori fiber space or a log minimal model $(Y,B_Y,\Mm)/Z$.
    \item If $K_Y+B_Y+\Mm_Y$ is nef$/Z$, then it is semi-ample$/Z$.
    \item If $(X,B,\Mm)$ is $\Qq$-factorial gdlt, then any $(K_X+B+\Mm_X)$-MMP$/Z$ with scaling of an ample$/Z$ $\Rr$-divisor terminates.
\end{enumerate}
\end{thm}
Theorem \ref{thm: gmm exists for g-crepant log structure} generalizes \cite[Theorem 1.1]{Bir12}(see also \cite[Theorem 1.6]{HX13}, \cite[Theorem 1.1]{Has19}) to the category of g-pairs. We remark that the authors proved Theorem \ref{thm: gmm exists for g-crepant log structure}(1)(3) in \cite[Theorem 1.3]{LX22} while Theorem \ref{thm: gmm exists for g-crepant log structure} completes the missing part (2). Finding this last missing piece is very important, as it allows us to deduce the existence of flips for glc g-pairs in full generality.

\begin{thm}\label{thm: glc flip exists}
Let $(X,B,\Mm)/U$ be an NQC glc g-pair and $f:X \to Z$ is a $(K_X+B+\Mm_X)$-flipping contraction$/U$. Then the flip $X^+\to Z$ of $f$ exists.
\end{thm}

Theorem \ref{thm: glc flip exists} removes the $\mathbb R$-Cartier condition of $\Mm_X$ as in \cite[Theorem 1.2]{HL21a}, hence gives a complete solution of \cite[Conjecture 3.12]{HL22}. We remark that the proof of Theorem \ref{thm: glc flip exists} is quite different from the proof of \cite[Theorem 1.2]{HL21a}. Indeed, the proof of Theorem \ref{thm: glc flip exists} gives an alternative proof of \cite[Theorem 1.2]{HL21a}.

\smallskip

The next result is the g-pair version of  \cite[Theorem 1.1]{HX13}(see also \cite[Theorem 1.4]{Bir12}, \cite[Theorem 1.1]{Has19}) in full generality. 

\begin{thm}\label{thm: gmm over U0 implies gmm over U}
Let $(X,B,\Mm)/U$ be an NQC glc g-pair and $U^0\subset U$ a non-empty open subset. Let $X^0:=X\times_UU^0$, $B^0:=B\times_UU^0$, and $\Mm^0:=\Mm\times_UU^0$. Assume that
\begin{enumerate}
\item $(X^0,B^0,\Mm^0)/U^0$ has a good minimal model, and
\item any glc center of $(X,B,\Mm)$ intersects $X^0$.
\end{enumerate}
Then $(X,B,\Mm)/U$ has a good minimal model.
\end{thm}

\begin{rem}
When $\Mm=\bm{0}$, Theorem \ref{thm: gmm over U0 implies gmm over U} is closely related to the properness of the moduli functor of stable schemes. Unfortunately, it seems difficult for us to apply Theorem \ref{thm: gmm over U0 implies gmm over U} in a similar way in the study of the moduli of g-pairs. In general, it is not clear if we can extend a glc structure on $X^0$ over $U^0$ to a glc structure on a compactification $X$ of $X^0$ over a compactification $U$ of $U^0$. This is mainly because a nef$/U^0$ divisor on $X^0$ usually does not extend to a nef divisor$/U$ on $X$. In fact, many properties for pairs in families do not hold for g-pairs anymore, see \cite{BH22} for examples where the theory of g-pairs presents extreme complications. We refer the reader to \cite{Bir22} for some new techniques about moduli for generalized pairs. 
\end{rem}

We remark that \cite[Theorem 1.1]{HL21a} proves Theorem \ref{thm: gmm over U0 implies gmm over U} under the additional assumption that $\Mm^0_{X^0}\sim_{\mathbb R,U^0}0$. The proof of Theorem \ref{thm: gmm over U0 implies gmm over U} is quite different from the proof of \cite[Theorem 1.1]{HL21a} as well. Indeed, the proof of Theorem \ref{thm: gmm over U0 implies gmm over U} also provides an alternative proof of \cite[Theorem 1.1]{HL21a}.

\medskip

The following result, which generalizes \cite[Theorem 1.5]{Bir12} to the category of g-pairs, is also important to the proofs of Theorems \ref{thm: gmm exists for g-crepant log structure}, \ref{thm: glc flip exists}, and \ref{thm: gmm over U0 implies gmm over U}. It is interesting to see that, although the finite generation of the generalized log canonical ring usually fails, it is still useful in the minimal model program for generalized pairs.

\begin{thm}\label{thm: finite generation imply semi-ample}
Let $(X,B,\Mm)/U$ be a $\Qq$-factorial gdlt $\Qq$-g-pair such that $f: X\to U$ is surjective. Let $U^0$ be a non-empty open set of $U$ and $X^0:=X\times_U U^0$. Assume that 
\begin{enumerate}
    \item $R(X/U,K_X+B+\Mm_X)$ is a finitely
generated $\Oo_U$-algebra, and
\item $(K_X+B+\Mm_X)|_{X^0}$ is semi-ample over $U^0$.
\end{enumerate}
Then $(X,B,\Mm)/U$ has a good minimal model. Moreover, any sequence of steps for the $(K_X+B+\Mm_X)$-MMP$/U$ with scaling of an ample$/U$ $\Rr$-divisor terminates with a good minimal model of $(X,B,\Mm)/U$.
\end{thm}

The key idea in the proofs of Theorems \ref{thm: gmm exists for g-crepant log structure}, \ref{thm: glc flip exists}, and \ref{thm: gmm over U0 implies gmm over U} is a Koll\'ar-type gluing theory which we will establish in Section \ref{sec: gluing}, combined with the minimal model program for special g-pairs as in \cite{LX22} (see also \cite{Has22}). As an important application of independent interest, we show that glc singularities are Du Bois. This is a generalization of \cite[Theorem 1.4]{KK10} to the category of generalized pairs, and will allow us to construct many Du Bois singularities without any log canonical structure (cf. Example \ref{ex: glc not lc}).

\begin{thm}\label{thm: glc sings are Du Bois}
Let $(X,B,\Mm)$ be an NQC glc g-pair. Then any union of glc centers of $(X,B,\Mm)$ is Du Bois. In particular, $X$ is Du Bois.
\end{thm}

We remark that \cite[Theorem 1.1]{FL22} shows that qlc (quasi-log canonical) singularities are Du Bois. Since any qlc pair is always a glc g-pair \cite[Remark 1.9]{Fuj22}, Theorem \ref{thm: glc sings are Du Bois} implies  \cite[Theorem 1.1]{FL22}.

\medskip

We expect the theorems above to have important applications in future studies of g-pairs. We state a few of them here. The first one is the extractability of non-canonical places of glc g-pairs:

\begin{thm}\label{thm: extracting divisor over glc structure}
Let $(X,B,\Mm)$ be an NQC glc g-pair, and $E$ a prime divisor that is exceptional over $X$ such that $a(E,X,B,\Mm)\in [0,1)$. Then there exists a birational morphism $f: Z\to X$ which extracts $E$ such that $-E$ is ample over $X$.
\end{thm}

When $\Mm=\bm{0}$, Theorem \ref{thm: extracting divisor over glc structure} is proved in \cite[Theorem 1]{Mor20}.

\medskip

We show the finite generation of the ring for any integral divisor which avoids glc centers:

\begin{thm}\label{thm: finite generation of R(X,A)}
Let $(X,B,\Mm)$ be an NQC glc g-pair, and $D$ an integral divisor on $X$, such that $\Supp D$ does not contain any glc center of $(X,B,\Mm)$. Then $R(X,D)$ is a finitely generated $\Oo_X$-algebra.
\end{thm}

When proving the theorems above, we get some counter-examples to some expected properties of g-pairs. We will summarize them in Section \ref{sec: ideas}. We hope they will be useful in future studies of generalized pairs. 

\medskip

\noindent\textit{Structure of the paper}. In Section \ref{sec: ideas}, we summarize our ideas of the proofs of the main theorems and provide some examples of g-pairs satisfying special properties. In Section \ref{sec: preliminaries}, we introduce some preliminary results that will be used in the rest of the paper. In Section \ref{sec: gluing}, we introduce the concept of \emph{glc crepant log structures}, a generalized pair version of crepant log structures for lc pairs \cite[Definition 4.28]{Kol13}, and establish a Koll\'ar-type gluing theory for this structure. In Section \ref{sec: key theorem}, we prove the key theorem Theorem \ref{thm: semi-ample over U0 implies semi-ample over U}. In Section \ref{sec: DB singularity} we explore the Du Bois property coming from glc crepant log structures and prove Theorem \ref{thm: of glc origin implies DB}. In Section \ref{sec: proof of the main theorems}, we use Theorem \ref{thm: semi-ample over U0 implies semi-ample over U} and Theorem \ref{thm: of glc origin implies DB} to prove our main theorems.

\medskip

\noindent\textbf{Acknowledgement}. The authors would like to thank their advisor Christopher D. Hacon for useful discussions and constant support. They would like to thank Jingjun Han, Yuchen Liu, and Chenyang Xu for useful discussions. We thank the referee for detailed suggestions. The second author is partially supported by NSF research grants no: DMS-1801851, DMS-1952522 and by a grant from the Simons Foundation; Award Number: 256202.

\section{Idea of the proof of Theorem \ref{thm: gmm exists for g-crepant log structure} and some examples}\label{sec: ideas}

It is clear that the main difficulty of the proof of Theorem \ref{thm: gmm exists for g-crepant log structure} will appear when gluing glc centers. For the usual pair case, there are two methods to resolve this issue:
\begin{enumerate}
    \item Show the fact that nef and log abundant implies semi-ample (cf. \cite{FG14,HX16,Has20a}).
    \item Show the finiteness of {\rm{\textbf{B}}}-representations (cf. \cite{FG14,HX16}).
\end{enumerate}
However, the nuances of glc g-pairs seem to pose some serious difficulties. Indeed, we will show later that both (1) and (2) have counter-examples.

\subsection{Idea of the proof of Theorem \ref{thm: gmm exists for g-crepant log structure}}

The key idea in our proof of Theorem \ref{thm: gmm exists for g-crepant log structure} is that, instead pursuing a more general statement as in the proofs of the usual pair case (like the finiteness of {\rm{\textbf{B}}}-representations), we shall fully utilize all conditions imposed and prove the finiteness of relations and the existence of geometric quotients only in this restricted setting. 

To better illustrate our idea, let's start with some cases when we can easily prove the finiteness of relations so that a ``direct" proof of gluing is possible. For example, suppose that $W$ is sdlt, $\pi: W^n\to W$ is the normalization or $W$, and $D^n$ is the double locus. Let $L_W$ be a semi-ample line bundle on $W$ which defines a contraction $g: W\to Y$. Let $g^n: W^n\to Y^n$ be the contraction induced by $L=\pi^*L_W$ and $T^n\rightrightarrows Y^n$ be the relation induced by the relation $D^n\rightrightarrows W^n$. Then the relation generated by $T^n\rightrightarrows Y^n$ is automatically finite, and the geometric quotient is just $Y$. 

This observation seems useless, as our goal --- the semi-ampleness of $L_W$, where $W$ is the non-gklt locus and $L$ is the restricted generalized log canonical divisor --- is already in the assumptions. Nevertheless, by applying induction on dimension, we may assume that $L_{W^n}:=\pi^*L_W$ is semi-ample. We have a key observation here: by a lemma of Koll\'ar  \cite[Lemma 9.55]{Kol13}, to prove the finiteness of relations, we only need to show the semi-ampleness of $L_W$ over a ``good" open subset of $W$. For arbitrary sdlt varieties $W$, or even if $W$ is the non-gklt locus of an arbitrary gdlt g-pair, such ``good" open subset may not exist. However, such good open set will automatically exist under the setting of Theorem \ref{thm: gmm over U0 implies gmm over U}, where we can let that open subset be the inverse image of $U^0$. 

Now the last thing we need to do is to establish a Koll\'ar-type gluing theory under the setting of Theorem \ref{thm: gmm over U0 implies gmm over U}. This is also not trivial: when a similar kind of Koll\'ar-type gluing theory was introduced in \cite{HX13,HX16} in the proof of the existence of lc flips, they ended up using the finiteness of \textbf{B}-representations which we want to avoid. Nevertheless, thanks to the MMPs developed in \cite{LX22} (see also \cite{HL21a,Has22}), we are able to reduce Theorem \ref{thm: gmm over U0 implies gmm over U} to the case when $(X^0,B^0,\Mm^0)/U^0$ is a good minimal model of itself (cf. Theorem \ref{thm: semi-ample over U0 implies semi-ample over U}). By the generalized canonical bundle formula  \cite{Fil20,FS20,HL21b,JLX22} and induction on dimension, we reduce to the case when $K_X+B+\Mm_X$ is big and nef (Step 2 of Theorem \ref{thm: semi-ample over U0 implies semi-ample over U}). Now we can get a gluing theory that can be directly applied for this case without using the finiteness of \textbf{B}-representations. More precisely, with the help of the generalized canonical bundle formula and the structure of $\mathbb P^1$-links for glc g-pairs \cite{FS20}, we may apply similar arguments as in \cite[Chapter 4]{Kol13} to establish a gluing theory for g-pairs with gdlt crepant log structures (see Section \ref{sec: gluing} for details). This eventually provides the gluing theory that we need, and all the main theorems will follow.

\subsection{Example}

In this subsection, we will provide three examples corresponding to three failed approaches towards Theorem \ref{thm: gmm exists for g-crepant log structure}. These approaches are:
\begin{enumerate}
    \item Try to get an lc structure on $X$ and show that a glc flip is also an lc flip.
    \item Try to show that nef and log abundant implies semi-ample (for g-pairs).
    \item Try to prove the finiteness of {\rm{\textbf{B}}}-representations (for g-pairs).
\end{enumerate}
All these three approaches are natural approaches and have played crucial roles before. Indeed, (2) and (3) are essentially used when proving the existence of lc flips \cite{Bir12,HX13}, while (1) is essentially used when proving the existence of $\Qq$-factorial glc flips \cite{HL21a}. We hope that our examples will illustrate some of the subtleties of working with glc g-pairs and be beneficial for future works.

\subsubsection{Glc pair without lc structure} A key observation in \cite{HL21a} indicates that, for any glc g-pair $(X,B,\Mm)/U$ such that $\Mm_X$ is $\mathbb R$-Cartier, any twist of $(X,B,\Mm)/U$ with any ample$/U$ $\mathbb R$-divisor will induce an lc structure on $X$ (cf. \cite[Lemma 5.18]{HL21a}). Actually, this observation leads to the proof of Theorem \ref{thm: glc flip exists} when $\Mm_X$ is $\mathbb R$-Cartier \cite[Theorem 1.2]{HL21a}. 

However, when dealing with non-$\mathbb Q$-factorial glc g-pairs, one cannot expect the existence of an lc structure on $X$ due to the following example:

\begin{ex}\label{ex: glc not lc}
Let $S$ be a projective lc variety such that $-K_S$ is nef but not big, and $(S,\Delta)$ is not lc for any $\Delta\in |-K_S|_{\mathbb R}$, i.e. $(S,0)$ does not have an $\mathbb R$-complement. Such $S$ exists, even if we additionally require that $S$ is smooth (cf. \cite[1.1 Example]{Sho00}, where $S=\PP_E(V)$ is a ruled surface over an elliptic curve $E$ and $V$ is a non-splitting vector bundle over $E$ of rank $2$). 

Let $L$ be an ample line bundle on $S$. Then the affine cone $Y:=C(S,L)$ is not potentially lc, i.e. for any $B_Y\geq 0$ on $Y$, $(Y,B_Y)$ is not lc. To see this, let $p:X:=BC(S,L)\to C(S,L)=Y$ be the blow-up of the vertex of $Y$ with exceptional divisor $E\simeq S$, then $\pi:BC(S,L)\to S$ is total space of the line bundle $L^{-1}$ over $S$ and $E$ is the zero section. If there exists $B_Y\ge0$ such that $(Y,B_Y)$ is lc, then
$$
p^*(K_Y+B_Y)=K_X+(1-a)E+B_X
$$
where $a\geq 0$ and $B_X:=f^{-1}_*B_Y$. Since $K_S$ is $\mathbb Q$-Cartier, $K_X$ is $\Qq$-Cartier. Since $\pi$ is smooth, we have $(K_X+E)|_E\sim_{\Qq} K_S$, hence
$$
K_S\sim_{\Rr}-(B_X|_{E}+aL).
$$
Since $-K_S$ is not big, $a=0$. In this case, $-K_S\sim_{\Rr}B_X|_{E}\ge0$ and $(S,B_X|_E)$ is lc by adjunction. This contradicts our assumption that  $(S,\Delta)$ is not lc for any $\Delta\in |-K_S|_{\mathbb R}$.

On the other hand, $Y$ does have a glc structure $(Y,0,\overline{M})$, where $M=\pi^*(-K_S)$ is a nef $\mathbb Q$-divisor on $X$. By Theorem \ref{thm: glc sings are Du Bois}, $Y$ is also an example of a variety which is Du Bois but has no lc structure.
\end{ex}

\subsubsection{Nef and log abundant do not imply semi-ample for g-pairs} By adopting and further developing the ideas of Hashizume \cite{Has20b,Has22}, in \cite{LX22}, we are able to prove Theorem \ref{thm: gmm exists for g-crepant log structure}(1)(3). We use the additional structure given by the morphism $f: X\rightarrow Z$ as in Theorem \ref{thm: gmm exists for g-crepant log structure}. In fact, by induction on dimension, we can reduce to the case when $(X,B,\Mm)$ is log abundant$/Z$. In the classical minimal model program, nefness and log abundance usually imply semi-ampleness (cf. \cite{FG14,HX16,Has20a}). Nef and abundant also imply semi-ample for gklt g-pairs: the first known proof of this result is \cite[Lemma 3.10]{Has22}; see also \cite[Theorem 2]{Cha22}.

However, we have the following example of a glc g-pair with nef and log abundant but not semi-ample generalized log canonical divisor:

\begin{ex}[{\cite[Example 1.4]{LX22}}]\label{ex: log abundant not semi-ample}
Let $C_0$ be a nodal cubic in $\Pp^2$ and $l$ the hyperplane class on $\Pp^2$. Let $P_1,P_2,...,P_{12}$ be twelve distinct points on $C_0$ which are different from the nodal point of $C_0$. Let
$$
\mu:X=\text{Bl}_{\{P_1,...,P_{12}\}}\to\Pp^2
$$
be the blow-up of $\Pp^2$ at the chosen points with the exceptional divisor $E=\sum_{i=1}^{12}E_i$, where $E_i$ is the prime exceptional divisor over $P_i$ for each $i$. Let $H:=\mu^*l$ and $C:=\mu^{-1}_*C_0$. Then $C\cong C_0$, $C\in|3H-E|$, and $K_X+C=\mu^*(K_{\Pp^2}+C_0)=0$. 

We consider the big divisor $M=4H-E\sim H+C$. Since $H$ is semi-ample and $M\cdot C=0$, $M$ is nef. Notice that $\Oo_C(M)=\Oo_{C_0}(4l-\sum_{i=1}^{12}P_i)$ and $\Pic^0(C)\cong\mathbb G_m$, where $\mathbb G_m$ is the multiplication group of $\Cc^*$. Let $\epsilon$ be any sufficiently small rational number, then $M-\epsilon C\sim_{\Qq}H+(1-\epsilon)C$ is ample by the Nakai-Moishezon Criterion. 

Suppose that $P_1,...,P_{12}$ are in general position such that $\Oo_C(M)$ is a non-torsion in $\Pic^0(C)$. Then $M$ is not semi-ample since $M|_C$ is not. However, the normalization $C^n$ of $C$ is $\Pp^1$, so $M|_{C^n}$ is semi-ample. We let $\Mm:=\overline{M}$ be the closure of $M$, i.e. $\Mm$ is the $\bb$-divisor such that $\Mm$ descends to $X$ and $\Mm_X=M$ (cf. \cite[Definition 2.9]{HL21a}). Then we have a glc g-pair $(X,C,\Mm:=\overline{M})$ such that both $M$ and $K_X+C+M$ are nef and log abundant with respect to $(X,C,\Mm)$, but $K_X+C+M$ is not semi-ample. 

Let $f: Y\to X$ be the blow-up at the node of $C_0$. Then $K_Y+C_1+C_2=f^*(K_X+C)$, where $C_2\cong \Pp^1$ is the $f$-exceptional divisor and $C_1\cong \Pp^1$ is the birational transform of $C$. We have that
\begin{enumerate}
    \item $(Y,C_1+C_2,\Mm)$ is a smooth gdlt g-pair,
    \item $K_Y+C_1+C_2+\Mm_Y=\Mm_Y=f^*M$ is nef and log abundant with respect to $(Y,C_1+C_2,\Mm)$,
    \item $(K_Y+C_1+C_2+\Mm_Y)|_{C_i}$ is semi-ample, and
    \item $K_Y+(1-\epsilon)C_1+(1-2\epsilon)C_2+\Mm_Y\sim f^*(M-\epsilon C)$ is big and semi-ample.
\end{enumerate}
However, $K_Y+C_1+C_2+\Mm_Y=f^*M$ is not semi-ample.
\end{ex}
 We also remark that conditions (1--4) in Example \ref{ex: log abundant not semi-ample} show that we will not be able to get any similar statement as \cite[Theorem 1.7]{Bir12}, \cite[Corollary 1.5]{HX16} for glc g-pairs, while those results are crucial in the proof of the existence of lc flips.

\subsubsection{\Bb-representations for g-pairs are not finite.} The main issue to prove Theorem \ref{thm: gmm exists for g-crepant log structure} is to glue the semi-ample structures on the glc centers together. For log canonical pairs, such gluing theory is established in \cite{FG14,HX16} thanks to the finiteness of \Bb-representations. Therefore, we want to investigate the finiteness of \Bb-representations for glc g-pairs as well. \cite{Hu21} indicated that the finiteness of \Bb-representations is expected to hold for g-pairs under some additional technical assumptions. However, we easily get the following very simple counter-example on the finiteness of \textbf{B}-representations for g-pairs.

\begin{ex}\label{ex: fail finiteness b representation}
Let $n$ be a positive integer and $(\Pp^n,0,\overline{M})$ a g-pair, where $M=(n+2)H\sim\Oo_{\Pp^n}(n+2)$ and $H$ is a hyperplane section on $\mathbb P^n$. Then the automorphisms of $\Pp^n$ which fix $H$ form an infinite subgroup $\Aut(\Pp^n,H)$ of $\Bir(\Pp^n,0,\overline{M})$. Since the representation of $\Aut(\Pp^n)\cong PGL(n+1,\Cc)$ on $H^0(\Pp^n, K_{\Pp^n}+M)\cong H^0(\Pp^n,\Oo_{\Pp^n}(1))$ is faithful, $\rho_1(\Bir(\Pp^n,0,\overline{M}))$ is infinite, where $\rho_m: \Bir(\Pp^n,0,\overline{M})\to \Aut(H^0(\Pp^n,mK_{\Pp^n}+mM))$.
\end{ex}

As a consequence of the failure of the finiteness of \textbf{B}-representations, the gluing theory for g-pairs is problematic. As shown in Example \ref{ex: relation is not finite in general} below, the semi-ampleness of a g-sdlt pair (cf. \cite{Hu21}) is quite subtle and is hard to distinguish from its normalization without any extra conditions.

\section{Preliminaries}\label{sec: preliminaries}

We adopt the standard notation and definitions in \cite{KM98,BCHM10} and will freely use them. 

\begin{defn}
Let $X\rightarrow U$ be a projective morphism and $D$ a Weil divisor on $X$ such that $|D/U|\not=\emptyset$. We let 
$$\Fix(D/U):=\sum_P(\inf_{D'\in |D/U|}\mult_PD')P$$ be the \emph{fixed part} of $D$, and let $\Mov(D):=D-\Fix(D)$ be the \emph{movable part} of $D$.
\end{defn}

\begin{defn}[Generalized pairs]\label{defn: gpairs}
For g-pairs, we adopt the same notation as in \cite{HL21a}. In particular, a generalized pair $(X,B,\Mm)/U$ consists of a normal quasi-projective variety $X$ associated with a projective morphism $X\rightarrow U$,
an $\Rr$-divisor $B$ on $X$, and a nef$/U$ $\bb$-divisor $\Mm$ over $X$, such that $K_X+B+\Mm_X$ is $\Rr$-Cartier. We make the following minor changes:
\begin{enumerate}
\item (Trivial glc centers) For any g-(sub-)pair $(X,B,\Mm)$, we will consider $X$ itself as a glc center and a non-gklt center of $(X,B,\Mm)$. $X$ will be called the \emph{trivial} glc center/\emph{trivial} non-gklt center of $(X,B,\Mm)$. We will let $\Ngklt(X,B,\Mm)$ be the union of all non-trivial non-gklt center of $(X,B,\Mm)$.
\item (Scheme structure of glc locus) We will always consider $\Ngklt(X,B,\Mm)$ as a scheme which is associated with the natural reduced scheme structure. In particular, if $(X,B,\Mm)$ is gdlt, then $\lfloor B\rfloor=\Ngklt(X,B,\Mm)$ is considered as both a divisor and a reduced scheme.
\item (Gplt) We say that a glc g-pair $(X,B,\Mm)/U$ is \emph{generalized plt} (\emph{gplt} for short) if $(X,B,\Mm)$ is gdlt and $\lfloor B\rfloor$ is normal.
\end{enumerate}
We also remark that different definitions of gdlt g-pairs in literature are now equivalent to each other thanks to \cite[Theorem 6.1]{Has22}.
\end{defn}

\begin{defn}\label{defn: b log abundance}
Let $(X,B,\Mm)/U$ be a sub-glc g-sub-pair and $D$ an $\Rr$-divisor on $X$. We say that $D$ is \emph{abundant}$/U$ if $\kappa_{\iota}(X/U,D)=\kappa_{\sigma}(X/U,D)$. We say that $D$ is \emph{log abundant}$/U$ with respect to $(X,B,\Mm)$ if $D$ is log abundant$/U$, and for any glc center $W$ of $(X,B,\Mm)$ with normalization $W^{\nu}$, $D|_{W^{\nu}}$ is abundant$/U$. We say that $(X,B,\Mm)$ is \emph{log abundant}$/U$ if $K_X+B+\Mm_X$ is log abundant$/U$ with respect to $(X,B,\Mm)$.
\end{defn}

\subsection{Perturbations of generalized pairs}

\begin{lem}[{cf. \cite[Proof of Lemma 4.4]{BZ16}, \cite[Page 717, Line 5]{HL22}}]\label{lem: gklt g-pair to pair}
Let $(X,B,\Mm)/U$ be a gklt g-pair and $f: Y\rightarrow X$ a birational morphism such that $\Mm$ descends to $Y$ and $\Mm_Y$ is big$/U$. Then there exists a klt pair $(X,\Delta)$ such that $K_X+B+\Mm_X\sim_{\mathbb R,U}K_X+\Delta$.
\end{lem}
\begin{proof}
Let $K_Y+B_Y+\Mm_Y:=f^*(K_X+B+\Mm_X)$. For any positive integer $n$, We may write $\Mm_Y=H_n+\frac{1}{n}E$ where $H_n$ is ample$/U$ and $E\geq 0$. Fix $n\gg 0$, then we may pick $A_n\in |H_n/U|_{\mathbb R}$ such that $(Y,B_Y+\frac{1}{n}E+A_n)$ is sub-gklt. We may let $\Delta:=f_*(B_Y+\frac{1}{n}E+A_n)$. 
\end{proof}

\begin{lem}\label{lem: reduction to Nlc locus}
Let $(X,B,\Mm)/U$ be a $\Qq$-factorial NQC gdlt g-pair. Assume that
\begin{enumerate}
    \item $L:=K_X+B+\Mm_X$ is nef$/U$ and big$/U$, 
    \item $W:=\Ngklt(X,B,\Mm)$, and
    \item $L|_W$ is semi-ample over $U$.
\end{enumerate}
Then $L$ is semi-ample over $U$.
\end{lem}
\begin{proof}
By the theory of Shokurov-type rational polytopes \cite[Theorem 2.28]{HL21a} (see also \cite[Proposition 3.16]{HL22}, \cite[Lemma 5.3]{HLS19},~\cite[Theorem 1.4]{Che20}) for generalized pairs, there exist real numbers $a_1,\dots,a_k\in (0,1]$ and $\Qq$-g-pairs $\{(X,B_i,\Mm^i)\}_{i=1}^k$ satisfying the following:
\begin{itemize}
    \item $\sum_{i=1}^ka_i=1$.
    \item $B=\sum_{i=1}^ka_iB_i$ and $\Mm=\sum_{i=1}^ka_i\Mm^i$.
    \item $(X,B_i,\Mm^i)$ is a gdlt $\Qq$-g-pair for any $i$.
    \item $L_i=K_X+B_i+\Mm_X^i$ is nef$/U$ and big$/U$.
    \item $L_i|_W$ is semi-ample$/U$.
    \item $\Ngklt(X,B_i,\Mm^i)=\Ngklt(X,B,\Mm)=W$ for each $i$.
\end{itemize}
Thus we may assume that $(X,B,\Mm)$ is a $\mathbb Q$-g-pair. (To see this, note that $(X,B_i,\Mm^i)/U,L_i$, and $W_i$ satisfying the conditions of Lemma \ref{lem: reduction to Nlc locus} and $L=\sum_{i=1}^ka_iL_i$. So $L$ is semi-ample over $U$ when $L_i$ is semi-ample for each $i$.)

Let $f:Y\to X$ be a log resolution of $(X,\Supp B)$ such that $\Mm$ descends to $Y$, and let $K_Y+B_Y+\Mm_Y:=f^*(K_X+B+\Mm_X)$. Since $L$ is nef$/U$ and big$/U$, we may write $L\sim_{\mathbb Q,U}H_n+\frac{1}{n}F$ for any positive integer $n$, such that $H_n\geq 0$ is ample and $F\geq 0$. Now for each $n$ and any positive integer $m$, we may write 
$$\Mm_Y+\frac{1}{2}f^*H_n\sim_{\mathbb Q,U}A_{n,m}+\frac{1}{m}E_{n},$$
where $A_{n,m}$ are ample$/U$ $\mathbb Q$-divisors and $E_n\geq 0$. For any $m\gg n\gg 0$, we have
$$\Nklt(Y,B_Y+\frac{1}{m}E_n+\frac{1}{n}f^*F)=\Nklt(Y,B_Y)=\Ngklt(Y,B_Y,\Mm),$$
thus we may pick $A_{n,m}\geq 0$ such that $$\Nlc(Y,\Delta_Y:=B_Y+A_{n,m}+\frac{1}{m}E_n+\frac{1}{n}f^*F)=\Ngklt(Y,B_Y,\Mm),$$
where $\Nlc(Y,\Delta_Y)$ is defined as in \cite[Section 7]{Fuj17}.
Let $\Delta:=f_*\Delta_Y$, then $\Delta\geq 0$, $\Nlc(X,\Delta)=\Ngklt(X,B,\Mm)=W$, and $2L-(K_X+\Delta)\sim_{\mathbb Q,U}\frac{1}{2}H_n$ is ample$/U$. The lemma follows from \cite[Theorems 4.5.5, 6.5.1]{Fuj17}, \cite[Theorem 5.3]{Amb03}.
\end{proof}

\begin{rem}\label{rem: to q coefficients}
As in the proof of Lemma \ref{lem: reduction to Nlc locus}, we will frequently use Shokurov-type rational polytopes to reduce g-pair questions to $\mathbb Q$-g-pair questions. To avoid redundancy, in the following, we will just cite \cite[Theorem 2.28]{HL21a} and do not list out all the details of the decomposition (e.g. we will not list out items from ``$\sum_{i=1}^ka_i=1$" to ``$\Ngklt(X,B_i,\Mm^i)=\Ngklt(X,B,\Mm)=W$ for each $i$" as in the proof of Lemma \ref{lem: reduction to Nlc locus}).
\end{rem}

The following result is an easy consequence of \cite[Lemma 5.18]{HL21a} although it is not in literature, so we write it here. We do not need it in the rest of the paper.

\begin{thm}\label{thm: bpf with mx r cartier}
Let $(X,B,\Mm)/U$ be a glc g-pair and $L$ a nef$/U$ Cartier divisor on $X$ such that $L-(K_X+B+\Mm_X)$ is ample$/U$. Assume that $\Mm_X$ is $\mathbb R$-Cartier. Then $mL$ is base-point-free$/U$ for any integer $m\gg 0$.
\end{thm}
\begin{proof}
Possibly replacing $\Mm$ with $(1-\epsilon)\Mm$ for some $0<\epsilon\ll 1$, we may assume that $\Ngklt(X,B,\Mm)=\Nklt(X,B)$. Let $A:=L-(K_X+B+\Mm_X)$. By \cite[Lemma 5.18]{HL21a}, there exists a birational morphism $h: W\rightarrow X$ such that $\Mm$ descends to $W$ and $\Supp(h^*\Mm_X-\Mm_W)=\Exc(h)$. We let $E:=h^*\Mm_X-\Mm_W$, then $E\geq 0$ and $E$ is $h$-exceptional.

Let $K_W+B_W:=h^*(K_X+B)$. By our construction, $\Exc(h)= \Supp E$ does not contain any lc place of $(X,B)$. Thus we may pick $E'\geq 0$ on $Y$ such that $-E'$ is ample$/X$ and $E'$ does not contain any lc place of $(X,B)$. Since $\Ngklt(X,B,\Mm)=\Nklt(X,B)$, we may find $0<\delta\ll 1$ such that $\frac{1}{2}h^*A-\delta E'$ is ample$/U$ and $(W,B_W+\delta E')$ is sub-lc. In particular, we may find an ample$/U$ $\Rr$-divisor $$0\leq H_W\sim_{\Rr,U}\Mm_W+\frac{1}{2}h^*A-\delta E'$$ on $W$ such that $(W,B_W+H_W+\delta E')$ is sub-lc. Let $\Delta:=B+h_*H_W$, then $(X,\Delta)$ is lc and $\Delta\sim_{\Rr,U}B+\Mm_X+\frac{1}{2}A$. In particular, $L-(K_X+\Delta)\sim_{\mathbb R,U}\frac{1}{2}A$ is ample$/U$. Theorem \ref{thm: bpf with mx r cartier} follows from \cite[Theorem 5.3]{Amb03}, \cite[Theorems 4.5.5, 6.5.1]{Fuj17}.
\end{proof}

\subsection{Canonical bundle formula}

We will follow the notation as in \cite{JLX22}. See also \cite{Fil20,FS20,HL21b} for related results.

\begin{defn}\label{defn contraction}
A \emph{contraction} is a projective morphism $f: Y\rightarrow X$ such that $f_*\mathcal{O}_Y=\mathcal{O}_X$. In particular, $f$ is surjective and has connected fibers.
\end{defn}

\begin{defn}[Glc-trivial fibration, {cf. \cite[Definition 2.10]{JLX22}}]
Let $(X,B,\Mm)/U$ be a g-sub-pair and $f: X\rightarrow Z$ a contraction$/U$. If
\begin{enumerate}
    \item $(X,B,\Mm)$ is sub-glc over the generic point of $Z$,
    \item $\rk f_*\mathcal{O}_X(\lceil\Aa^*(X,B,\Mm)\rceil)=1$, and
    \item $K_X+B+\Mm_X\sim_{\Rr,Z}0$,
\end{enumerate}
then we say that $f: (X,B,\Mm)\rightarrow Z$ is a glc-trivial fibration$/U$. 
\end{defn}

\begin{defn}
Let $(X,B,\Mm)/U$ be an NQC g-sub-pair and $f: (X,B,\Mm)\rightarrow Z$ is a glc-trivial fibration$/U$, and $(Z,B_Z,\NN)$ an NQC g-sub-pair on $Z$. We say that $(Z,B_Z,\NN)$ is an (NQC) g-sub-pair \emph{induced by the canonical bundle formula$/U$ of} $f: (X,B,\Mm)\rightarrow Z$ if $K_X+B+\Mm_X\sim_{\mathbb R}f^*(K_Z+B_Z+\NN_Z)$ and $a(D,Z,B_Z,\NN)=1-t_D(X,B,\Mm;f)$ for any prime divisor $D$ over $Z$, where $t_D(X,B,\Mm;f)$ are glc thresholds defined as in \cite[Definition 2.12]{JLX22}. 

By \cite[Theorem 2.23]{JLX22}, if $B\geq 0$ over the generic fiber of $f$, then there always exists an NQC g-sub-pair induced by the canonical bundle formula$/U$ of $f: (X,B,\Mm)\rightarrow Z$. Moreover, it is not hard to see that if $(X,B,\Mm)$ is a $\Qq$-g-sub-pair, then the induced g-sub-pair on $Z$ can also be chosen as a $\Qq$-g-sub-pair. We will frequently use these facts in the rest of the paper.
\end{defn}

\subsection{Crepant log structures}

\begin{defn}
A \emph{glc crepant log structure} is of the form $f: (X,B,\Mm)\rightarrow Z$, where
\begin{enumerate}
    \item $(X,B,\Mm)/Z$ is a glc g-pair,
    \item $K_X+B+\Mm_X\sim_{\Rr,Z}0$, and
    \item $f$ is a contraction. In particular, $f_*\Oo_X=\Oo_Z$.
\end{enumerate}
In addition, if
\begin{enumerate}
    \item[(4)] $(X,B,\Mm)$ is gdlt, 
\end{enumerate}
then we say that $f: (X,B,\Mm)\rightarrow Z$ is a \emph{gdlt crepant log structure}. An \emph{NQC glc (resp. gdlt) crepant log structure} is a glc (resp. gdlt) crepant log structure $f: (X,B,\Mm)\rightarrow Z$ such that $\Mm$ is NQC$/Z$.

We remark that glc crepant log structures are also known as \emph{generalized log Calabi-Yau fibrations}. We use the wording ``glc crepant log structure" because we mainly use this structure for Koll\'ar's glueing theory (see Section \ref{sec: gluing}) while \cite[Definition 4.28]{Kol13} uses the wording ``crepant log structure".

For any irreducible subvariety $W\subset Z$, we say that $W$ is a \emph{glc center} of a glc crepant log structure $f: (X,B,\Mm)\rightarrow Z$, if there exists a glc center $W_X$ of $(X,B,\Mm)$ such that $W=f(W_X)$. For any (not necessarily closed) point $z\in Z$, we say that $z$ is a \emph{glc center} of $f: (X,B,\Mm)\rightarrow Z$ if $\bar z$ is a glc center of $f: (X,B,\Mm)\rightarrow Z$.
\end{defn}

\begin{lem}\label{lem: glc centers come from cbf}
Let $(X,B,\Mm)/U$ be an NQC glc g-pair, $f: (X,B,\Mm)\rightarrow Z$ a glc-trivial fibration$/U$, and $(Z,B_Z,\NN)/U$ an NQC g-pair induced by the canonical bundle formula of $f: (X,B,\Mm)\rightarrow Z$. Then for any irreducible subvariety $W$ of $Z$, $W$ is a glc center of $f: (X,B,\Mm)\rightarrow Z$ if and only if $W$ is a glc center of $(Z,B_Z,\NN)$.
\end{lem}
\begin{proof}
The if part follows \cite[Theorem 2.23]{JLX22}  and the only if part follows from \cite[Theorem 2.16(2)]{LX22}.
\end{proof}

\begin{defn}
Let $(X,B,\Mm)$ and $(X',B',\Mm')$ be two g-pairs. We say that  $(X,B,\Mm)$ and $(X',B',\Mm')$ are \emph{crepant equivalent to each other} if there exist birational morphisms $p: W\rightarrow X$ and $q: W\rightarrow X'$ such that $\Mm'=\Mm$ and $p^*(K_X+B+\Mm_X)=q^*(K_{X'}+B'+\Mm'_{X'})$.
\end{defn}

\subsection{\texorpdfstring{$\mathbb P^1$}{}-links}
We recall the definition and results on $\mathbb P^1$-links as in \cite{FS20}. This is a generalization of \cite[Theorem 4.40]{Kol13} to the category of generalized pairs. We partially refine the definitions (e.g. we define $\mathbb P^1$-links for $\mathbb R$-g-pairs) to make our arguments more clear and general.

\begin{defn}[Standard $\mathbb P^1$-link, cf. {\cite[Definition 2.21]{FS20}}]\label{defn: standard p1 link}
A \emph{standard $\mathbb P^1$-link} is a glc g-pair $(X,B,\Mm)/Z$ with a projective morphism $f:X\to T$ over $Z$ satisfying the following properties.
\begin{enumerate}
\item $K_X+B+\Mm_X\sim_{\mathbb R,T}0$,
\item there exists a birational morphism $X'\rightarrow X$ such that $\Mm_{X'}\sim_{\mathbb R,T}0$,
\item $\lfloor B\rfloor=D_1+D_2$, where $D_1,D_2$ are prime divisors and $f|_{D_i}: D_i\rightarrow T$ are isomorphisms,
\item $(X,B,\Mm)/Z$ is gplt, and
\item every reduced fiber of $f$ is isomorphic to $\mathbb P^1$.
\end{enumerate}
We call $D_1$ and $D_2$ the \emph{horizontal sections} of $(X,B,\Mm)/T$.
\end{defn}

\begin{defn}[$\mathbb P^1$-link, cf. {\cite[Definition 2.23]{FS20}}]\label{defn: p1 link}
Let $(X,B,\Mm)/Z$ be a gdlt g-pair associated with a projective morphism $f: X\rightarrow Z$, such that $K_X+B+\Mm_X\sim_{\mathbb R,Z}0$. Let $Z_1$, $Z_2$ be two glc centers of $(X,B,\Mm)$. We say that $Z_1$ and $Z_2$ are \emph{directly $\mathbb P^1$-linked$/Z$} if there exists an irreducible subvariety $W\subset X$, such that either $W$ is a glc center of $(X,B,\Mm)$ or $W=X$, and we have the following. Let $(W,B_W,\Mm^W)/Z$ be a gdlt g-pair induced by adjunction to the higher-codimensional glc center $W$, i.e.
$$K_W+B_W+\Mm^W_W:=(K_X+B+\Mm_X)|_W,$$
such that
\begin{enumerate}
    \item $Z_i\subset W$ for each $i$,
    \item $f(W)=f(Z_1)=f(Z_2)$, and
    \item there exists a g-pair $(W',B_{W'},\Mm^W)$ crepant equivalent to $(W,B_W,\Mm^W)$ and a projective morphism $h: W'\rightarrow T$ over $Z$, such that $(W',B_{W'},\Mm^W)/T$ is a $\mathbb P^1$-link and $Z_1|_{W'},Z_2|_{W'}$ are the horizontal sections of $(W',B_{W'},\Mm^W)/T$.
\end{enumerate}
We say that $Z_1$ and $Z_2$ are \emph{$\mathbb P^1$-linked$/Z$} if either $Z_1=Z_2$, or there exists an integer $n\geq 2$ and glc centers $Z_1',\dots,Z_n'$ of $(X,B,\Mm)$, such that $Z_1'=Z_1,Z_n'=Z_2$, and $Z'_i$ and $Z'_{i+1}$ are directly $\mathbb P^1$-linked$/Z$ for any $1\leq i\leq n-1$.
\end{defn}

\begin{thm}[cf. {\cite[Theorem 3.5]{Bir20}, \cite[Theorem 1.4]{FS20}}]\label{thm: P1 link for gdlt crepant log structure} 
Let $(X,B,\Mm)/U$ be an NQC gdlt g-pair associated with a projective morphism $f: X\rightarrow U$, such that $K_X+B+\Mm_X\sim_{\mathbb R,U}0$. Let $s\in U$ be a (not necessarily closed) point such that $f^{-1}(s)$ is connected (as a $k(s)$-scheme). Let
$$\mathcal{S}:=\{V\mid V\text{ is a glc center of }(X,B,\Mm), s\in f(V)\}$$
and $Z,W\in\mathcal{S}$ be two elements such that $Z$ is minimal in $\mathcal{S}$ with respect to the inclusion. Then there exists $Z_W\in\mathcal{S}$ such that $Z_W\subset W$, and $Z$ and $Z_W$ are $\mathbb P^1$-linked$/U$. In particular, any minimal elements in $\mathcal{S}$ with respect to inclusion are $\mathbb P^1$-linked$/U$ to each other.
\end{thm}
\begin{proof}
It following from \cite[Theorem 2.28]{HL21a} and \cite[Theorem 1.4]{FS20}.
\end{proof}

\begin{lem}\label{lem: glc locus is unibranch}
Let $f: (X,B,\Mm)\rightarrow Z$ be an NQC glc crepant log structure and $z\in Z$ a (not necessarily closed) point. Let $$\mathcal{S}_z:=\{V\mid V\text{ is a glc center of }f: (X,B,\Mm)\rightarrow Z, z\in V\}.$$
Then:
\begin{enumerate}
    \item There exists a unique element $W\in\mathcal{S}_z$ that is minimal with respect to inclusion. 
    \item $W$ is unibranch at $z$, i.e.  the completion $\hat{W}_z$ is irreducible.
    \item Any intersection of glc centers of $f: (X,B,\Mm)\rightarrow Z$  is also a union of glc centers.
\end{enumerate}
\end{lem}
\begin{proof}
The proof is exactly the same as in \cite[Proof of Corollary 4.41]{Kol13} except that we replace \cite[Theorem 4.40]{Kol13} with Theorem \ref{thm: P1 link for gdlt crepant log structure}. For the reader's convenience, we give a full proof here. 

Possibly replacing $(X,B,\Mm)$ with a gdlt model, we may assume that $(X,B,\Mm)$ is gdlt. For any any element $W\in\mathcal{S}_z$ that is minimal, there exists a glc center $Z_W$ of $(X,B,\Mm)$ that is minimal among all glc centers whose image on $Z$ is equal to $W$ with respect to inclusion. By Theorem \ref{thm: P1 link for gdlt crepant log structure}, all such $Z_W$ are $\mathbb P^1$-linked$/Z$ to each other, hence their images on $Z$ are the same. This proves (1). (2) follows from (1) by considering every \'etale neighborhood of $z$. 

For any glc centers $W_1,W_2$ on $Z$, let $z\in W_1\cap W_2$ be any point, and $W$ the unique element minimal element of $\mathcal{S}_z$. Then $W\subset W_1\cap W_2$, and we get (3).
\end{proof}

The following lemma should be well-known, but we cannot find any reference.
\begin{lem}\label{lem: inversion of adjunction}
Let $(X,B,\Mm)/Z$ be an NQC gdlt g-pair, $S$ a component of $\lfloor B\rfloor$, and $(S,B_S,\Mm^S)/Z$ the gdlt g-pair induced by the adjunction 
$$K_S+B_S+\Mm^S_S:=(K_X+B+\Mm_X)|_S.$$
Then:
\begin{enumerate}
    \item Any glc center of $(S,B_S,\Mm^S)$ is a glc center of $(X,B,\Mm)$.
    \item Any glc center of $(X,B,\Mm)$ that is contained in $S$ is a glc center of $(S,B_S,\Mm^S)$.
\end{enumerate}
\end{lem}
\begin{proof}
By \cite[Theorem 6.1]{Has22}, there exists a log resolution $f: \tilde X\rightarrow X$ of $(X,\Supp B)$ and an open subset $X^0\subset X$, such that $\Mm$ descends to $\tilde X$, $X^0$ contains the generic point of any glc center of $(X,B,\Mm)$, and $f$ is an isomorphism over $X^0$. Let $K_{\tilde X}+\tilde B+\Mm_{\tilde X}:=f^*(K_X+B+\Mm_X)$ and let $\tilde S$ be the strict transform of $S$ on $\tilde X$, then
$f|_{\tilde S}$ is a log resolution of $(S,\Supp B_S)$ such that $\Mm^S$ descends to $\tilde S$, i.e.
$$f|_{\tilde S}^*(K_S+B_S+\Mm^S_S)=K_{\tilde S}+B_{\tilde S}+\Mm^S_{\tilde S}:=(K_{\tilde X}+\tilde B+\Mm_{\tilde X})|_{\tilde S}.$$
Thus any glc center of $(S,B_S,\Mm^S)$ is a glc center of $(\tilde S,B_{\tilde S},\Mm^S)$, hence a glc center of $(\tilde X,\tilde B,\Mm)$, and hence a glc center of $(X,B,\Mm)$, which shows (1). On the other hand, any glc center of $(X,B,\Mm)$ that is contained in $S$ is a glc center of $(\tilde X,\tilde B,\Mm)$ that is contained in $\tilde S$, hence a glc center of $(\tilde S,B_{\tilde S},\Mm^S)$, and hence a glc center of $(S,B_S,\Mm^S)$, which shows (2).
\end{proof}

\begin{lem}\label{lem: gdlt crepant log structure is compatible under subadjunction}
Let $f: (X,B,\Mm)\to Z$ be an NQC gdlt crepant log structure and $Y\subset X$ a glc center. Let
$$
f|_Y: Y\xrightarrow{f_Y}Z_Y\xrightarrow{\pi} Z
$$
be the Stein factorization of $f|_Y$, and $(Y,B_Y,\Mm^Y)/Z$ the NQC gdlt g-pair induced by adjunction to the higher-codimensional glc center $Y$, i.e. 
$$K_Y+B_Y+\Mm_Y^Y:=(K_X+B+\Mm_X)|_Y.$$
Then:
\begin{enumerate}
\item $f_Y: (Y,B_Y,\Mm^Y)\rightarrow Z_Y$ is a gdlt crepant log structure.
\item For any glc center $W_Y\subset Z_Y$ of $f_Y: (Y,B_Y,\Mm^Y)\rightarrow Z_Y$, $\pi(W_Y)$ is a glc center of $f: (X,B,\Mm)\rightarrow Z$.
\item For any glc center $W\subset Z$ of $f: (X,B,\Mm)\rightarrow Z$, every irreducible component of $\pi^{-1}(W)$ is a glc center of  $f_Y: (Y,B_Y,\Mm^Y)\rightarrow Z_Y$.
\end{enumerate}
\end{lem}

\begin{proof}
The proof is exactly the same as in \cite[Corollary 4.42]{Kol13} except that we use Theorem \ref{thm: P1 link for gdlt crepant log structure} in replace of \cite[Theorem 4.40]{Kol13}. We also have a proof of (3) in \cite[Proof of 4.1]{JLX22}. For the reader's convenience, we give a full proof here.

(1) We only need to show that $(Y,B_Y,\Mm^Y)$ is gdlt, which follows from \cite[Lemma 2.6]{HL22}. 

(2) There exists a glc center $V_Y$ of $(Y,B_Y,\Mm^Y)$ such that $f_Y(V_Y)=W_Y$. By Lemma \ref{lem: inversion of adjunction}, $V_Y$ is also a glc center of $(X,B,\Mm)$. Thus $\pi(W_Y)=f(V_Y)$ is a glc center of $f: (X,B,\Mm)\rightarrow Z$.

(3) Let $z$ be the generic point of $W$. Since the question is \'etale local, possibly replacing $Z$ by an \'etale neighborhood of $z$ and replacing $Y$ with its irreducible components, we may assume that $f^{-1}(z)\cap Y$ is connnected, and we only need to show that there exists a glc center $V_Y$ of $f_Y: (Y,B_Y,\Mm^Y)\rightarrow Z_Y$ such that $f_Y(V_Y)$ is an irreducible component of $\pi^{-1}(W)$.

Let $V_X$ be a minimal glc center of $(X,B,\Mm)$ which dominates $W$, i.e. $V_X$ is minimal in
$$\{V\mid V\text{ is a glc center of }(X,B,\Mm), V\text{ dominates }W\}$$
with respect to inclusion. Then $f(V_X)=W$. By Theorem \ref{thm: P1 link for gdlt crepant log structure}, there exists a glc center $V_Y\subset Y$ of $(X,B,\Mm)$ that is $\mathbb P^1$-linked$/Z$ to $V_X$. By Lemma \ref{lem: inversion of adjunction}, $V_Y$ is also a glc center of $(Y,B_Y,\Mm^Y)$. Thus $f_Y(V_Y)\subset Z_Y$ is a glc center of $f_Y: (Y,B_Y,\Mm^Y)\rightarrow Z_Y$. Moreover, since $V_Y$ is $\mathbb P^1$-linked$/Z$ to $V_X$, $f(V_Y)=f(V_X)=W$. Thus $f_Y(V_Y)$ is an irreducible component of $\pi^{-1}(W)$ and we are done.
\end{proof}

\begin{rem}
In the setting of Lemma \ref{lem: gdlt crepant log structure is compatible under subadjunction}, $f_Y$ actually induces an NQC glc structure $(Z_Y,B_{Z_Y},\Mm^{Z_Y})$ on $Z_Y$ by the canonical bundle formula, and also induces an NQC glc structure $(T,B_{T},\Mm^{T})$ on the normalization $T$ of $f(Y)$ by \cite[Theorem 1.2]{HL21b}. Let $(Z,B_Z,\Mm^Z)$ be an NQC glc g-pair induced by the canonical bundle formula$/Z$ of $f: (X,B,\Mm)\rightarrow Z$, then we can also perform sub-adjunction by \cite[Theorem 5.1]{HL21b} to $T$, and the induced structure will coincide with $(T, B_{T},\Mm^{T})$ up to an $\mathbb R$-linear equivalence of the moduli part (see \cite[Section 4]{JLX22}).
\end{rem}

\section{Koll\'ar-type gluing theory for generalized pairs}\label{sec: gluing}

In this section we will review Koll\'ar’s powerful gluing theory of finite quotients. We refer for \cite[Section 5 and Section 9]{Kol13} for more details. We will develop the gluing theory we need for glc crepant log structures in this section.

In this section, we will generally choose the notation $(X,\Delta,\Mm)$ instead of $(X,B,\Mm)$ for g-pairs, as $B$ is used in the boundary of stratifications.

\subsection{Definitions}

\begin{defn}[{\cite[Definition 9.15]{Kol13}}] 
Let $X$ be a scheme. A {\it stratification} of $X$ is a decomposition of $X$ into a finite disjoint union of reduced locally closed subschemes. We will consider stratifications where the strata are of pure dimensions
and are indexed by their dimensions. We write $X=\cup_{i}S_iX$ where $S_iX\subset X$ is the $i$-th
dimensional stratum. Such a stratified scheme is denoted by $(X,S_*)$. We also
assume that $\cup_{i\le j}S_iX$ is closed for every $j$. The {\it boundary} of $(X,S_*)$ is the closed subscheme
$$
B(X,S_*):=\cup_{i<\dim X}S_iX=X\backslash S_{\dim X}X,
$$
and is denoted by $B(X)$ if the stratification $S_*$ is clear. 

Let $(X, S_*)$ and $(Y, S_*)$ be stratified schemes. We say that $f:X\to Y$ is a {\it stratified morphism} if $f(S_iX)\subset S_iY$ for every $i$. Since $S_iX$ are disjoint with each other, $f: X\to Y$ is a stratified morphism if and only if $S_iX=f^{-1}(S_iY)$.

Let $(Y, S_*)$ be a stratified scheme and $f:X\to Y$ a quasi-finite morphism such that $f^{-1} (S_iY)$ has pure dimension $i$ for every $i$ . Then $S_iX:=f^{-1}(S_iY)$ defines a stratification of $X$. We denote it by $(X,f^{-1}S_*)$, and we say that $f:X\to(Y,S_*)$ is \emph{stratifiable}.
\end{defn}

\begin{defn}[{\cite[Definition 9.16]{Kol13}}]
Let $(X, S_*)$ be stratified variety. A relation $(\sigma_1,\sigma_2): R\rightrightarrows (X,S_*)$ is {\it stratified} if each $\sigma_i$ is stratifiable and $\sigma_1^{-1}S_*=\sigma_2^{-1}S_*$. Equivalently,
there exists a stratification $(R,\sigma^{-1}S_i)$, such that $r\in\sigma^{-1}S_iR$ if and only if $\sigma_1(r)\in S_iX$ and if and only if $\sigma_2(r)\in S_iX$.
\end{defn}

\begin{defn}[{\cite[Definition 9.18]{Kol13}}]
Let $(X,S_*)$ be a stratified scheme such that $X$ is an excellent scheme. The normality conditions (N), (SN), (HN), and (HSN)  are defined in the following ways.
\begin{enumerate}
    \item[(N)] We say that $(X,S_*)$ has {\it normal strata}, or that it satisfies (N), if each $S_iX$ is normal.
    \item[(SN)] We say that $(X,S_*)$ has {\it semi-normal boundary}, or that it satisfies (SN), if $X$ and $B(X,S_*)$ are both semi-normal.
    \item[(HN)] We say that $(X,S_*)$ has {\it hereditarily normal strata}, or that it satisfies (HN), if \begin{enumerate}
            \item the normalization $\pi: (X^n,\pi^{-1}S_*)\to (X,S_*)$ is stratifiable,
            \item $(X^n,S_*^n)$ satisfies (N), and
            \item $B(X^n,\pi^{-1}S_*)$ satisfies (HN).
    \end{enumerate}       
    \item[(HSN)] We say that $(X,S_*)$ has {\it hereditarily semi-normal boundary}, or that it
    satisfies (HSN), if \begin{enumerate}
            \item the normalization $\pi: (X^n,\pi^{-1}S_*)\to (X,S_*)$ is stratifiable,
            \item $(X,S_*)$ satisfies (SN), and
            \item $B(X^n,\pi^{-1}S_*)$ satisfies (HSN).
    \end{enumerate}
\end{enumerate}
\end{defn}

Next we give a special stratification that is induced by the glc crepant log structure. 

\begin{defn}[Glc stratification]
Let $f:(X,\Delta,\Mm)\to Z$ be a glc crepant log structure. Let $S^*_i(Z,X,\Delta,\Mm)\subset Z$ be the union of all $\le i$-dimensional glc centers of $f:(X,\Delta,\Mm)\to Z$, and
$$
S_i(Z,X,\Delta,\Mm):=S^*_i(Z,X,\Delta,\Mm)~\backslash ~S^*_{i-1}(Z,X,\Delta,\Mm).
$$
If the glc crepant log structure $f:(X,\Delta,\Mm)\to Z$ is clear from the context, we will use $S_i(Z)$ for abbreviation. It is clear that each $S_i(Z)$ is a locally closed subspace of $Z$ of pure dimension $i$, and $Z$ is the disjoint union of all $S_i(Z)$. 

The stratification of $Z$ induced by $S_i(Z)$ is called the \emph{generalized log canonical stratification} (\emph{glc stratification} for short) of $Z$ induced by $f:(X,\Delta,\Mm)\to Z$.  Since this is the only stratification we are going to use in the rest of this paper, we usually will not emphasize the glc crepant structure $f:(X,\Delta,\Mm)\to Z$, and we will denote the corresponding stratified scheme by $(Z,S_*)$. The \emph{boundary} of $(Z,S_*)$ is the closed subspace
$$B(Z,S_*):=Z\backslash S_{\dim Z}(Z)=\cup_{i<\dim Z}S_i(Z).$$
\end{defn}

\begin{defn}\label{defn: of glc origin}
We say that a semi-normal stratified space $(Y,S_*)$ is \textit{of generalized log canonical (glc) origin} if $S_i(Y)$ is unibranch for any $i$, and there are glc crepant log structures $f_j:(X_j,\Delta_j,\Mm^j)\to Z_j$ with glc stratifications $(Z_j,S_{*}^j)$ and a finite surjective stratified morphism $\pi: \amalg_j(Z_j,S_{*}^j)\to (Y,S_*)$. Moreover, if $f_j:(X_j,\Delta_j,\Mm^j)\to Z_j$ are NQC glc repant log structures, then we say that $(Y,S_*)$ is \textit{of NQC glc origin}.
\end{defn}

\subsection{Basic properties}

The following theorem and its proof are very similar to \cite[Proposition 4.32]{Kol13}.
\begin{thm}\label{thm: glc locus is semi-normal}
Let $f:(X,\Delta,\Mm)\to Z$ be an NQC glc crepant log structure. Let $W\subset Z$ be the union of all glc centers of $f:(X,\Delta,\Mm)\to Z$ except $Z$, and $B(W)\subset W$ the union of all non-maximal (with respect to inclusion) glc centers that are contained in $W$. Then
\begin{enumerate}
    \item $W$ is semi-normal, and
    \item $W\backslash B(W)$ is normal.
\end{enumerate}
\end{thm}
\begin{proof}
By \cite[Theorem 2.28]{HL21a}, we may assume that $(X,\Delta,\Mm)$ is a $\Qq$-g-pair. Let $(Z,\Delta_Z,\NN)/U$ be a glc $\mathbb Q$-g-pair induced by the canonical bundle formula$/U$ of $f: (X,\Delta,\Mm)\rightarrow Z$. By Lemma \ref{lem: glc centers come from cbf}, the glc centers of $(Z,\Delta_Z,\NN)$ are exactly the glc centers of $f: (X,\Delta,\Mm)\rightarrow Z$. Possibly replacing $(X,\Delta,\Mm)$ with a gdlt model of $(Z,\Delta_Z,\NN)$, we may assume that $f$ is birational and $(X,\Delta,\Mm)$ is $\mathbb Q$-factorial gdlt. We have $W=f(\lf\Delta\rf)$. Let $\Delta':=\{\Delta\}$. We consider the exact sequence 
$$
0\to\Oo_X(-\lf\Delta\rf)\to\Oo_X\to\Oo_{\lf\Delta\rf}
$$
and its push-forward
$$
\Oo_Z=f_*\Oo_X\to f_*\Oo_{\lf\Delta\rf}\stackrel{\delta}{\longrightarrow}R^1f_*\Oo_X(-\lf\Delta\rf).
$$
By Lemma \ref{lem: gklt g-pair to pair}, we can find a $\Qq$-divisor $\Delta''\ge 0$ such that $$-\lf\Delta\rf\sim_{\Qq,Z}K_X+\Delta'+\Mm_X\sim_{\Qq,Z}K_X+\Delta''$$ and $(X,\Delta'')$ is klt. By \cite[Corollary 10.40]{Kol13}, $R^if_*\Oo_{X}(-\lf\Delta\rf)$ is torsion free for every $i$. On the other hand, $f_*\Oo_{\lf\Delta\rf}$ is supported on $W$, hence it is a torsion sheaf. Thus the connecting map $\delta$ is zero, hence $\Oo_Z\twoheadrightarrow f_*\Oo_{\lf\Delta\rf}$ is surjective. Since this map factors through $\Oo_W$, we conclude that $\Oo_W\twoheadrightarrow f_*\Oo_{\lf\Delta\rf}$ is also surjective, hence an isomorphism.

Note that $\lf\Delta\rf$ has only nodes at codimension 1 points and it is $S_2$ by \cite[Corollary 2.88]{Kol13}. By \cite[Lemma 10.14]{Kol13}, $\lf\Delta\rf$ is semi-normal. By \cite[Lemma 10.15]{Kol13}, $W$ is semi-normal. This is (1).

To prove (2), let $V\subset\lf\Delta\rf$ be an irreducible component of its non-normal locus. Then $V$ is an lc center of $(X, \Delta)$, hence a glc center of $(X,\Delta,\Mm)$. Thus $f(V)\subset Z$ is a glc center. Hence either $f(V)$ is an irreducible component of $W$, or $f(V)\subset B(W)$. Thus \cite[Complement 10.15.1]{Kol13} implies that $W \backslash B(W)$ is normal.
\end{proof}

Theorem \ref{thm: glc locus is semi-normal} has the following interesting corollary. We do not need it in the rest of the paper.

\begin{cor}
Let $(X,\Delta,\Mm)$ be an NQC glc g-pair. Then $\Ngklt(X,\Delta,\Mm)$ is semi-normal.
\end{cor}
\begin{proof}
It follows from Theorem \ref{thm: glc locus is semi-normal} when $f$ is the identity morphism.
\end{proof}

\begin{lem}\label{lem: (Z,S) is U and SN}(cf. \cite[Lemma 5.26]{Kol13})
Let $f:(X,\Delta,\Mm)\to Z$ be an NQC glc crepant log structure and $(Z,S_*)$ the induced glc stratification. Then
\begin{enumerate}
    \item  $S_i(Z)$ is unibranch for every $i$, and 
    \item  $B(Z,S_*)$ is semi-normal.
\end{enumerate}
\end{lem}

\begin{proof}
(1) follows from Lemma \ref{lem: glc locus is unibranch}(2) and (2) follows from Theorem \ref{thm: glc locus is semi-normal}.
\end{proof}

\begin{lem}\label{lem: stratification is compatible under adjunction} (cf. \cite[Proposition 4.42]{Kol13})
Let $f: (X,\Delta,\Mm)\to Z$ be an NQC gdlt crepant log structure, $(Z,S_*)$ its induced glc stratification, and $Y\subset X$ a glc center of $(X,\Delta,\Mm)$. Let $(Y,\Delta,\Mm^Y)/Z$ be the NQC gdlt g-pair induced by adjunction to higher-codimensional glc center $Y$, i.e. 
$$K_Y+\Delta_Y+\Mm^Y_Y:=(K_X+\Delta+\Mm_X)|_Y,$$
We consider the Stein factorization of $f|_Y$
$$(Y,\Delta_Y,\Mm^Y)\stackrel{f_Y}{\longrightarrow}W\stackrel{\pi}{\longrightarrow}Z.$$
Then:
\begin{enumerate}
    \item $f_Y:(Y,\Delta_Y,\Mm^Y)\to W$ is an NQC gdlt crepant log structure which induces a glc stratification $(W,S_*)$.
    \item $S_i(W)=\pi^{-1}(S_i(Z))$ for every $i$.
\end{enumerate}
\end{lem}
\begin{proof}
It follows from Lemma \ref{lem: gdlt crepant log structure is compatible under subadjunction}.
\end{proof}

\begin{thm}\label{thm: (Z,S) is HN and HSN}
Let $f:(X,\Delta,\Mm)\to Z$ be an NQC glc crepant log structure and $(Z,S_*)$ the induced glc stratification. Then $(Z,S_*)$ satisfies (HN) and (HSN).
\end{thm}
\begin{proof}
By Lemma \ref{lem: (Z,S) is U and SN} and \cite[Definitions 9.18,~9.19]{Kol13}, $(Z,S_*)$ satisfies (HU) and (HSN). By \cite[Theorem 9.21]{Kol13}, $(Z,S_*)$  satisfies (HN).
\end{proof}

\begin{lem}(cf. \cite[5.29]{Kol13})\label{lem: glc stratification is of glc origin}
Every NQC glc stratification is of NQC glc origin. More precisely, let $f:(X,\Delta,\Mm)\to W $ be an NQC glc crepant log structure and $Y\subset W$ any union of glc centers. Then $(Y, S_*)$ is of NQC glc origin, where $S_i(Y)=Y\cap S_i(W)$ for each $i$.
\end{lem}

\begin{proof}
By Theorem \ref{thm: (Z,S) is HN and HSN} and \cite[Theorem 9.26]{Kol13} we know that $Y$ is semi-normal and $S_i(Y)$ is unibranch for each $i$. Then we can apply Lemma \ref{lem: stratification is compatible under adjunction} to each glc center of $f: (X,\Delta,\Mm)$ contained in $Y$ to conclude that $(YS_*)$ is of NQC glc origin.
\end{proof}

\subsection{Constructions of glc stratifications}

\begin{cons}[Gluing theory of glc crepant structures]\label{cons: gluing part 1}
Let $(X,\Delta,\Mm)/U$ be an NQC gdlt g-pair, $W\subset\lf\Delta\rf$ a reduced divisor, $\pi: W^n\rightarrow W$ the normalization of $W$, $D$ the double locus of $W^n$, $D^n$ the normalization of $D$, $\tau: D^n\rightarrow D^n$ the induced involution, and $(\tau_1,\tau_2): D^n\rightrightarrows W^n$ a finite stratified equivalence relation whose normalization map is given by the quotient morphism $\pi: W^n\to W=W^n/R$, where $R$ is the finite equivalence relation generated by $D^n$. 

Let $L_W:=(K_X+\Delta+\Mm_X)|_W$, $$L:=(K_X+\Delta+\Mm_X)|_{W^n}=K_{W^n}+\Delta_{W^n}+\Mm^{W^n}_{W_n}$$
where $(W^n,\Delta_{W^n},\Mm^{W^n})/U$ is the NQC gdlt g-pair by adjunction to $W^n$, and suppose that $L$ is semi-ample$/U$. Let $g^n: W^n\rightarrow Y^n$ and $h^n: D^n\rightarrow T^n$ be the morphisms$/U$ induced by $L$ and $L|_{D^n}$ respectively so that we have the commutative diagram
\begin{displaymath}
    \xymatrix{ 
        D^n \ar[dd]_{h^n}\ar@<.5ex>[rr]^{\tau_1} \ar@<-.5ex>[rr]_{\tau_2} && W^n \ar[dd]^{g^n} \\
        &&\\
        T^n \ar@<.5ex>[rr]^{\sigma_1}\ar@<-.5ex>[rr]_{\sigma_2} && Y^n 
    }
\end{displaymath}
where $(\sigma_1,\sigma_2): T^n\rightrightarrows Y^n$ are induced by $(\tau_1,\tau_2): D^n\rightrightarrows W^n$. We let $(D^n,\Delta_{D^n},\Mm^{D^n})/U$ be the gdlt g-pair induced by the adjunction
$$
K_{D^n}+\Delta_{D^n}+\Mm^{D^n}_{D^n}=(K_{W^n}+\Delta_{W^n}+\Mm^{W^n}_{W_n})|_{D^n}.
$$
It is clear that $g^n:(W^n,\Delta_{W^n},\Mm^{W^n})\to Y^n$ and  $h^n: (D^n,\Delta_{D^n},\Mm^{D^n})\rightarrow T^n$ are gdlt crepant log structures. We let $(Y^n,S_*(Y^n))$ and $(T^n,S_*(T^n))$ be their induced stratified schemes respectively.
\end{cons}

\begin{cons}\label{cons: glue part 2}
Notations and conditions as in Construction \ref{cons: gluing part 1}. Assume that $(X,B,\Mm)$ is a $\mathbb Q$-g-pair. Let $m$ be a sufficiently divisible positive integer such that $mL_W$ is Cartier, $|mL/U|$ defines $g^n$, and there exists a very ample$/U$ divisor $H$ on $Y^n$ such that $(g^{n})^*H=M$.

Let $p_W: W^n_M\to W^n$, $p_Y: Y^n_H\to Y^n$ be the total spaces of the line bundles $M$ and $H$ respectively. Let $\Delta_{W^n_M}:=p_W^{-1}(\Delta_{W^n})$, and $g^n_M: (W^n_M,\Delta_{W^n_M},p^*_W\Mm^{W^n})\to Y^n_H$ the gdlt crepant log structure with induced stratification $(Y^n_H, S_*(Y^n_H):=p_Y^{-1}S_*(Y^n))$. 

Let $p_D: D^n_M\to D^n$ and $p_T: T^n_H\to T^n$ be the total spaces of the line bundles $M|_{D^n}$ and $H|_{T^n}$. Let $\Delta_{D^n_M}:=p_D^{-1}(\Delta_{D^n})$, and $h^n_M: (D^n_M,\Delta_{D^n_M},p^*_D\Mm^{D^n})\to T^n_H$ the gdlt crepant log structure with induced stratification $(T^n_H, S_*(T^n_H):=p_T^{-1}S_*(T^n))$. 

Then we have a finite pre-relation $(\sigma_{1H},\sigma_{2H}): T^n_H\rightrightarrows Y^n_H$ induced by the finite relation $(\tau_{1M},\tau_{2M}): D^n_M\rightrightarrows W^n_M$, where $\tau_{1M},\tau_{2M}: D^n_M\rightarrow W^n_M$ are liftings of $\tau_{1},\tau_{2}$ respectively.
\end{cons}

\begin{lem}[{cf.~\cite[Lemma 3.11]{HX13}}]\label{lem: induced relation is stratified}
Notations and conditions as in Construction \ref{cons: gluing part 1}. Then
\begin{enumerate}
    \item $(\sigma_1,\sigma_2): T^n\rightrightarrows Y^n$ gives a stratified equivalence relation, and
    \item $(Y^n,S_*(Y^n))$ and $(T^n,S_*(T^n))$ satisfy (HN) and (HSN).
\end{enumerate}
If we have the additional notations and conditions as in Construction \ref{cons: glue part 2}, then\begin{enumerate}
    \item[(3)] $(\sigma_{1H},\sigma_{2H}): T^n_H\rightrightarrows Y^n_H$ gives a stratified equivalence relation, and
    \item[(4)] $(Y^n_H,S_*(Y^n_H))$ and $(T^n_H,S_*(T^n_H))$ satisfy (HN) and (HSN).
\end{enumerate}
\end{lem}
\begin{proof}
(2)(4) follow from Theorem \ref{thm: (Z,S) is HN and HSN}. We prove
(1)(3). For any glc center $V$ of $(D^n,\Delta_{D^n},\Mm^{D^n})$ (resp. of $(D^n_M,\Delta_{D^n_M},p^*_D\Mm^{D^n})$), $\tau(V)$ (resp. $\tau_M(V)$) is also a glc center on $D^n$ (resp. $D^n_M$). Thus the glc stratification induced by $h^n: (D^n,\Delta_{D^n},\Mm^{D^n})\to T^n$ (resp. $h^n_M: (D^n_M,\Delta_{D^n_M},p_D^*\Mm^{D^n})\to T^n_H$) is the same as the glc stratification induced by $h^n\circ\tau:(D^n,\Delta_{D^n},\Mm^{D^n})\to T^n$ (resp. $h^n_M\circ\tau_M: (D^n_M,\Delta_{D^n_M},p_D^*\Mm^{D^n})\to T^n_H$). Hence we only need to check that $\sigma^{-1}S_*(Y^n)$ (resp. $\sigma_H^{-1}S_*(Y^n_H)$) coincides with $S_*(T^n)$ (resp. $S_*(T^n_H)$), where $\sigma$ (resp. $\sigma_H$) is the canonical morphism $T^n\to Y^n$ (resp. $T^n_H\to Y^n_H$). But this follows directly from Lemma \ref{lem: stratification is compatible under adjunction}. 
\end{proof}

\subsection{Remarks and an example}
Notations and conditions as in Construction \ref{cons: glue part 2}. If $\Mm=0$, then $(W,\Delta_W)$ is sdlt. By \cite[Section 4]{HX16}, both $T^n\rightrightarrows Y^n$ and $T^n_H\rightrightarrows Y^n_H$ generate finite equivalence relations. By \cite[Theorem 9.21]{Kol13}, the geometric quotients $Y=Y^n/T^n$ and $Y_H=Y^n_H/T^n_H$ exist. Possibly by replacing $m$ with a multiple, $Y_H$ is a line bundle over $Y$, whose pullback to $W$ is exactly $mL_W$. In general, the pro-finite equivalence relation generated by $T^n\rightrightarrows Y^n$ and $T^n_H\rightrightarrows Y^n_H$ can be described as some almost group actions (\cite[Definition 9.32]{Kol13}) which is actually given by some crepant birational subgroup on the glc centers. Thanks to the finiteness of {\bf B}-representation for lc pairs \cite{HX13,FG14}, these groups are finite, hence the relations are also finite. 

However, when $\Mm\neq0$ and $(W,\Delta_W,\Mm^W)$ is only g-sdlt (cf. \cite{Hu21}), one should not expect that finiteness still holds without extra conditions or structures. We have already shown the failure of the finiteness of \textbf{B}-representations (cf. Example \ref{ex: fail finiteness b representation}). The following example will show that
\begin{enumerate}
    \item the relation generated by $T^n\rightrightarrows Y^n$ may not be finite and the geometric quotient $Y^n/T^n$ may not exist, and
    \item the relation generated by $T^n_H\rightrightarrows Y^n_H$ may not be finite, even when the geometric quotient $Y^n/T^n$ exists. 
\end{enumerate}

\begin{ex}\label{ex: relation is not finite in general}
(1) Let $\lambda\in\Cc^*$ and consider $\Pp^1\times \mathbb{A}^1$, which can be regarded as the total space of a trivial line bundle over $\Pp^1$. We define $\phi_\lambda: \{0\}\times\mathbb{A}^1\simeq\{\infty\}\times\mathbb{A}^1$ by $(0,t)\mapsto(\infty,\lambda t)$ and glue $\{0\}\times\mathbb{A}^1$ and $\{\infty\}\times\mathbb{A}^1$ together using $\phi_\lambda$ to get a demi-normal variety $M$ with projection $p:M\to C$, where $C$ is a nodal cubic. Then $M$ is a total space of a line bundle (also denoted by $M$) on $C$. Moreover, $M\in\Pic^0(C)\simeq\mathbb{G}_m=\Cc^*$ and can be canonically regarded as $\lambda\in\Cc^*$. Then:
\begin{itemize}
    \item $W:=C$ is sdlt and $K_C\sim 0$. We let $\Mm^W:=\overline{M}$.
    \item $\pi:\Pp^1\to C$ is the normalization and $D^n\rightrightarrows\Pp^1$ is the involution of two points $\{0,\infty\}$.
    \item $g^n: W^n\to Y^n$ is just $\Pp^1\to\Spec~\Cc$, and $T^n\rightrightarrows Y^n$ is trivial and finite. Therefore, the geometric quotient $Y^n/T^n$ exists and is equal to $\Spec~\Cc$. 
\end{itemize}
But from the line bundle aspect, we have the following:
\begin{itemize}
    \item $\pi_M:W^n_M\to M$ is $\Pp^1\times\mathbb{A}^1\to M$.
    \item $D^n\rightrightarrows\Pp^1\times\mathbb{A}^1$ is induced by $\phi_{\lambda}:\{0\}\times\mathbb{A}^1\simeq\{\infty\}\times\mathbb{A}^1$.
    \item $H$ is trivial and $g^n: W^n_M\to Y^n_H$ is the projection $\Pp^1\times\mathbb{A}^1\to \mathbb{A}^1$.
    \item $T^n\rightrightarrows \mathbb{A}^1$ is given by $\phi_\lambda,\phi^{-1}_\lambda$, and $\id$. Therefore, the relation generated by $T^n\rightrightarrows\mathbb{A}^1$ can viewed as the cyclic group $\langle\lambda\rangle\subset\Cc^*$, which is finite if and only if $\lambda$ is a root of unity.
\end{itemize}
(2) We can also compactify the above total spaces of line bundles to get projective examples when $Y^n/T^n$ does not exist.

Let $W:=\PP_C(\Oo_C\oplus M)$ be a $\Pp^1$-bundle over $C$, and let $C'\subset W$ be the section at infinity, which belongs to $|\Oo_W(1)|$. Then $W^n=\Pp^1\times\Pp^1$, and $g^n:W^n\to Y^n$ is the second projection $p_2: \Pp^1\times\Pp^1\to\Pp^1$.

Notice that $K_W$ is Cartier since $W$ is a locally complete intersection. Let $N:=3C'$, then $$\pi^*(K_W+N)=K_{W^n}+\{0\}\times\Pp^1+\{\infty\}\times\Pp^1+\pi^*N=p_2^*(\{\infty\})$$ 
is semi-ample. $N$ is nef since $\pi^*N$ is nef, so we see that $(W,0,\overline{N})$ is g-sdlt. However, the relation generated by $T^n\rightrightarrows\Pp^1$ is given by
$$\{[x,y]\sim[x',y']|[x',y']=[x,\lambda^ly] \text{ for some $l$}\}$$
and is finite if and only if $\lambda$ is a root of unity.
\end{ex}

\section{From gluing theory to abundance}\label{sec: key theorem}
The goal of this section is to prove the following theorem: 

\begin{thm}[{cf. \cite[Theorem 4.1]{HX13}}]\label{thm: semi-ample over U0 implies semi-ample over U}
Let $(X,B,\Mm)/U$ be a $\mathbb Q$-factorial NQC gdlt g-pair, $U^0$ a non-empty subset of $U$, and $X^0:=X\times_UU^0$. Assume that
\begin{enumerate}
    \item any glc center of $(X,B,\Mm)$ intersects $X^0$,
    \item $K_X+B+\Mm_X$ is nef$/U$, and
    \item $(K_X+B+\Mm_X)|_{X^0}$ is semi-ample$/U^0$.
\end{enumerate}
Then $K_X+B+\Mm_X$ is semi-ample$/U$. In particular, $(X,B,\Mm)/U$ is a good minimal model of itself.
\end{thm}

Before we prove Theorem \ref{thm: semi-ample over U0 implies semi-ample over U}, we need to prove Theorem \ref{thm: finite generation imply semi-ample}.

\begin{proof}[Proof of Theorem \ref{thm: finite generation imply semi-ample}]
Since termination and semi-ampleness$/U$ are both local on $U$, we can assume that $U$ is affine.

Let $m$ be a sufficiently divisible positive integer such that $m(K_X+B+\Mm_X)$ is Cartier and $m(K_X+B+\Mm_X)|_{X^0}$ is base-point-free$/U^0$, which defines a contraction$/U^0$ $h^0: X^0\to V^0$. Since $R(X/U,K_X+B+\Mm_X)$ is a finitely
generated $\Oo_U$-algebra, possibly replacing $m$ with a multiple, there
exist a log resolution $g: W\to X$ of $(X,\Supp B)$, a Weil divisor $E\ge 0$ on $W$, and a base-point-free$/U$ divisor $F$ on $W$, such that $\Mm$ descends to $W$,
$$
\Fix(g^*(lm(K_X+B+\Mm_X))/U)=lE, \text{and}  \Mov(g^*(lm(K_X+B+\Mm_X))/U)=lF
$$
for any positive integer $l$. Let $h: W\rightarrow V$ be the
contraction$/U$ defined by $|lF|$. Since $m(K_X+B+\Mm_X)|_{X^0}$ is base point free$/U^0$ and defines $h^0$, $V\times_UU^0=V^0$, and $E$ is vertical over $V$.

Let $B_W:=g^{-1}_*B+\Exc(g)_{\text{red}}$. Then $(W,B_W,\Mm)$ is a log smooth model of $(X,B,\Mm)$. We have
$$
m(K_W+B_W+\Mm_W)=g^*m(K_X+B+\Mm_X)+E'
$$
where $E'\ge0$ is exceptional over $X$. Thus
$$
\Fix(lm(K_W+B_W+\Mm_W)/U)=lE+lE', \text{ and }\Mov(lm(K_W+B_W+\Mm_W)/U)=lF.
$$

Let $B^0:=B\times_UU^0, B_W^0:=B_W\times_UU^0$, and $\Mm:=\Mm\times_UU^0$. We run a $(K_W+B_W+\Mm_W)$-MMP$/V$ with scaling of an ample divisor. Since $K_{X^0}+B^0+\Mm^0_{X^0}$ is semi-ample$/U^0$ and $K_{X^0}+B^0+\Mm^0_{X^0}\sim_{\mathbb Q,V^0}0$, $(X^0,B^0,\Mm^0)/U^0$ is a weak glc model of $(W^0,B^0_W,\Mm^0)/U^0$ and $(X^0,B^0,\Mm^0)/V^0$ is a weak glc model of $(W^0,B^0_W,\Mm^0)/V^0$. By \cite[Lemma 3.15]{HL21a}, $(W^0,B^0_W,\Mm^0)/V^0$ has a log minimal
model. By \cite[Theorem 2.24]{HL21a}, the $(K_W+B_W+\Mm_W)$-MMP$/V$ terminates over $V^0$. Let $\phi: W\dashrightarrow Y'$ be the induced birational map$/V$. 

Let $B_{Y'},E_{Y'},E'_{Y'}$, and $F_{Y'}$ be the strict transforms of $B_W,E,E'$, and $F$ on $Y'$ respectively. Since
$$ m(K_W+B_W+\Mm_W)\sim E+E'+F\sim E+E'$$
over $V^0$, $E_{Y'}+E'_{Y'}\sim_\Qq 0$ over $V^0$. In particular, $E_{Y'}+E'_{Y'}$ is vertical over $V$. 

Since $\phi: W\dashrightarrow Y'$ is a partial $(K_W+B_W+\Mm_W)$-MMP, 
$$\Fix((lE_{Y'}+lE'_{Y'}+lF_{Y'})/U)=\Fix(g^*(lm(K_{Y'}+B_{Y'}+\Mm_{Y'}))/U)=l(E_{Y'}+E_{Y'}').$$

By \cite[Lemma 3.2]{Bir12}, $E_{Y'}+E'_{Y'}$ is very exceptional (cf. \cite[Definition 3.1]{Bir12}) over $V$. By \cite[Proposition 3.8]{HL22}, we may run a $(K_{Y'}+B_{Y'}+\Mm_{Y'})$-MMP$/V$ with scaling of an ample divisor which terminates with a a log minimal model $(Y,B_Y,\Mm)/V$, such that
$$
m(K_Y+B_Y+\Mm_Y)\sim_{\Qq,V}E_Y+E'_Y=0,
$$
where $E_Y$ and $E_Y'$ are the strict transforms of $E_{Y'}$ and $E'_{Y'}$ on $Y$ respectively.
In particular, $m(K_Y+B_Y+\Mm_Y)\sim_{\Qq,U}F_Y$. Thus $K_Y+B_Y+\Mm_Y$ is semi-ample$/U$, hence $(Y,B_Y,\Mm_Y)/U$ is a good log minimal model of $(W,B_W,\Mm)/U$. By \cite[Lemma 3.10]{HL21a}, $(Y,B_Y,\Mm_Y)/U$ is a good log minimal model of $(X,B,\Mm)/U$. The moreover part of the theorem follows from \cite[Theorem 2.24, Lemma 3.9]{HL21a}.
\end{proof}

\begin{proof}[Proof of Theorem \ref{thm: semi-ample over U0 implies semi-ample over U}] Since semi-ampleness$/U$ is local on $U$, we can assume that $U$ is affine. By \cite[Theorem 2.28]{HL21a}, we may assume that $(X,B,\Mm)/U$ is a $\Qq$-g-pair. We let $B^0:=B\times_UU^0$ and $\Mm^0:=\Mm\times_UU^0$.

We may apply induction on dimensions. When $\dim X=1$ the theorem is obvious. Thus we may assume that $\dim X=d$ for some integer $d\geq 2$, and assume that the theorem holds in dimension $\leq d-1$. In particular, we may assume that $(K_X+B+\Mm_X)|_S$ is semi-ample$/U$ for any glc center $S$ of $(X,B,\Mm)$.

\medskip

\noindent\textbf{Step 1}. In this step, we construct an auxiliary g-pair $(V,B_V,\NN)/U$.

\smallskip

 Let $m>0$ be a sufficiently divisible integer such that $m(K_X+B+\Mm_X)$ is Cartier and $|m(K_X+B+\Mm_X)|_{X^0}|$ is base-point-free$/U^0$, which defines a contraction $h^0: X^0\to V^0$ over $U^0$. Let $h: X\dashrightarrow V$ be an Iitaka fibration$/U$ of $m(K_X+B+\Mm_X)$, then $h|_{X^0}=h^0$ is a morphism. We let $g: Y\to X$ be a log resolution of $(X,\Supp B)$ such that $\Mm$ descends to $Y$ and the induced birational map $Y\dashrightarrow V$ is a morphism. We can write $$K_Y+B_Y+\Mm_Y=g^*(K_X+B+\Mm_X)+E,$$
 where $B_Y\geq 0, E\geq 0$, and $B_Y\wedge E=0$. Then $E$ is exceptional over $X$, $(X,B,\Mm)/U$ is a weak glc model of $(Y,B_Y,\Mm)/U$, and the image of any glc center of $(Y,B_Y,\Mm)$ intersects $U^0$.
 
 By \cite[Theorem 4.2]{LX22} and \cite[Theorem 2.28]{HL21a}, we have the following commutative diagram
 \begin{center}$\xymatrix{
Y'\ar@{->}[d]_{h'}\ar@{.>}[r]^{f}& Y\ar@{->}[r]^{g}\ar@{->}[dr] & X\ar@{.>}[d]^{h}\\
V'\ar@{->}[rr]^{\varphi} & & V\\
}$
\end{center}
satisfying the following conditions:
\begin{itemize}
\item $h'$ is a contraction, $f$ is birational, and $\varphi: V'\rightarrow V$ is a resolution of $V$.
\item $(Y',B_{Y'},\Mm)$ is a $\mathbb Q$-factorial gdlt $\mathbb Q$-g-pair.
\item $K_{Y'}+B_{Y'}+\Mm_{Y'}\sim_{\mathbb Q,V'}0$.
\item Any weak glc model of $(Y,B_Y,\Mm)/U$ is a weak glc model of $(Y',B_{Y'},\Mm)/U$. In particular, $(X,B,\Mm)/U$ is a weak glc model of $(Y',B_{Y'},\Mm)/U$.
\item Any weak glc model of $(Y^0,B^0_Y,\Mm^0)/U$ is a weak glc model of $(Y'^0,B^0_{Y'},\Mm^0)/U$, where $Y^0:=Y\times_UU^0,Y'^0:=Y'\times_UU^0,B^0_Y:=B_Y\times_UU^0$, $B^0_{Y'}:=B_{Y'}\times_UU^0$, and $\Mm^0:=\Mm\times_UU^0$. In particular, $(X^0,B^0,\Mm^0)/U^0$ is a weak glc model of $(Y'^0,B^0_{Y'},\Mm^0)/U$.
\item Any glc center of $(Y',B_{Y'},\Mm)$ intersects  $Y'^0$.
\end{itemize}
By \cite[Theorem 2.16]{LX22}, there exists a glc $\mathbb Q$-g-pair $(V',B_{V'},\NN)/U$ induced by the canonical bundle formula$/U$ of $h': (Y',B_{Y'},\Mm)\rightarrow V'$, such that the image of any glc center of $(V',B_{V'},\NN)$ in $U$ intersects $U^0$. Since $h$ is an Iitaka fibration$/U$ of $K_X+B+\Mm_X$, $K_{V'}+B_{V'}+\NN_{V'}$ is big$/U$.

\medskip

\noindent\textbf{Step 2}. In this step, we reduce to the case when $K_X+B+\Mm_X$ is big$/U$.

\smallskip

We let $B^0:=B\times_UU^0$. Since $(X^0,B^0,\Mm^0)/U^0$ is a weak glc model of $(Y'^0,B_{Y'}^0,\Mm^0)/U^0$, there exist two birational morphisms $p: X''\rightarrow Y'$ and $q: X''\rightarrow X$, such that
$$p^*(K_{Y'}+B_{Y'}+\Mm_{Y'})|_{Y'^0}=q^*(K_X+B+\Mm_X)|_{X^0}+E^0$$
where $E^0\geq 0$ is exceptional over $X$ \cite[Lemma 3.8]{HL21a}.

By construction, $K_{X^0}+B^0+\Mm^0_{X^0}\sim_{\mathbb Q,V^0}0$. Since $K_{Y'}+B_{Y'}+\Mm_{Y'}\sim_{\mathbb Q,V'}0$, $K_{Y'^0}+B^0_{Y'}+\Mm^0_{Y'^0}\sim_{\mathbb Q,V'^0}0$. Since $V'\rightarrow V$ is birational, $V'\cong V$ over the generic point of $V$. Thus over the generic point of $V$,
$$K_{Y'}+B_{Y'}+\Mm_{Y'}\sim_{\mathbb Q}0\sim_{\mathbb Q}K_{X}+B+\Mm_{X},$$
and $(X,B,\Mm)$ is a good minimal model of $(Y',B_{Y'},\Mm)$. Thus $(X,B,\Mm)$ and $(Y',B_{Y'},\Mm)$ are crepant over the generic point of $V$. 

Since the $\Qq$-equivalence class of the moduli part of the canonical bundle formula only depends on the generic fiber of the fibration and canonical bundle formulas are compatible with base change, there exists a glc g-pair $(V^0,B^0_{V},\NN^0)/U^0$ induced by the canonical bundle formula of $h^0: X^0\rightarrow V^0$, such that $\NN^0=\NN\times_UU^0$. Let $V'^0:=V'\times_UU^0$ and $B^0_{V'}:=B_{V'}\times_UU^0$. Then
\begin{align*}
   &(h'|_{Y'^0})^*((K_{V'^0}+B^0_{V'}+\NN^0_{V'^0})-(\varphi|_{V'^0})^*(K_{V^0}+B^0_{V}+\NN^0_{V^0}))\\
    =&(K_{Y'^0}+B^0_{Y'}+\Mm^0_{Y'^0})-p_*q^*(K_X+B+\Mm_X)|_{X^0}=E^0\geq 0
\end{align*}
is exceptional over $X^0$. Therefore, $$0\leq (K_{V'^0}+B^0_{V'}+\NN^0_{V'^0})-(\varphi|_{V'^0})^*(K_{V^0}+B^0_{V}+\NN^0_{V^0})$$
is exceptional over $V^0$. Since $K_{V^0}+B^0_{V}+\NN^0_{V^0}$ is ample$/U^0$, $(V^0,B^0_{V},\NN^0)/U^0$ is a weak glc model of $(V'^0,B^0_{V'},\NN^0)/U^0$. By \cite[Lemmas 3.9, 3.15]{HL21a}, $(V'^0,B^0_{V'},\NN^0)/U^0$ has a good minimal model.

Let $(\tilde{V},B_{\tilde{V}},\NN)$ be a gdlt model of $(V',B_{V'},\NN)$, $\tilde{V}^0:=\tilde{V}\times_UU^0$, and $B^0_{\tilde V}:=B_{\tilde V}\times_UU^0$. By \cite[Theorem 3.14]{HL21a}, $(\tilde V^0,B^0_{\tilde V},\NN^0)/U^0$ has a good minimal model. By \cite[Lemma 2.7]{LX22} and \cite[Lemmas 3.9]{HL21a}, we may run a partial $(K_{\tilde V}+B_{\tilde V}+\NN_{\tilde V})$-MMP$/U$ $(\tilde V,B_{\tilde V},\NN)\dashrightarrow (\widehat V,B_{\widehat V},\NN)$, such that $(K_{\tilde V}+B_{\tilde V}+\NN_{\tilde V})|_{\widehat V^0}$ is semi-ample$/U^0$, where $\widehat V^0:=\widehat V\times_UU^0$. Now we run a $(K_{\widehat V}+B_{\widehat V}+\NN_{\widehat V})$-MMP$/U$ with scaling of an ample divisor
$$(\widehat V,B_{\widehat V},\NN)=(V_0,B_{V_0},\NN)\dashrightarrow (V_1,B_{V_1},\NN)\dashrightarrow\dots\dashrightarrow (V_i,B_{V_i},\NN)\dashrightarrow\dots.$$ 
Then the induced birational map $\widehat V\dashrightarrow V_i$ is an is an isomorphism over $U^0$. Since the image of any glc center of $(\widehat V,B_{\widehat V},\NN)$ on $U$ intersects $U^0$, the image of any glc center of $(V_i,B_{V_i},\NN)$ on $U$ intersects $U^0$. By induction hypothesis, $(V_i,B_{V_i},\NN)$ is log abundant$/U$ for each $i$. By \cite[Theorem 7.6]{LX22} (cf. \cite[Theorem 3.15]{Has22} when $X,U$ are projective varieties), this MMP terminates with a log minimal model $(\bar V,B_{\bar V},\NN)/U$ of $(V',B_{V'},\NN)/U$. Moreover, the image of any glc center of  $(\bar V,B_{\bar V},\NN)/U$ intersects $U^0$. By construction, 
\begin{align*}
    R(X/U,K_X+B+\Mm_X)&=R(Y'/U,K_{Y'}+B_{Y'}+\Mm_{Y'}) && \text{(Weak glc model)}\\
                    &=R(V'/U,K_{V'}+B_{V'}+\NN_{V'}) && \text{(Pullback)}\\
                  &=R(\bar{V}/U,K_{\bar{V}}+B_{\bar{V}}+\NN_{\bar{V}}) && (\text{Gdlt model+MMP}).
\end{align*}
If $\dim \bar V<\dim X$, then by induction hypothesis,  $K_{\bar{V}}+B_{\bar{V}}+\NN_{\bar{V}}$ is semi-ample$/U$, hence $R(\bar{V}/U,K_{\bar{V}}+B_{\bar{V}}+\NN_{\bar{V}})$ is finitely generated, so $R(X/U,K_X+B+\Mm_X)$ is finitely generated, and the theorem follows from Theorem \ref{thm: finite generation imply semi-ample}. Thus we may assume that $\dim \bar V=\dim X$, hence $K_X+B+\Mm_X$ is big$/U$.

\medskip

\noindent\textbf{Step 3}. We use gluing theory in Section \ref{sec: gluing} to prove the theorem. 

\smallskip

We let $W:=\lf B\rf=\Ngklt(X,B,\Mm)$, $W^0:=W\times_UU^0$, $L_W:=(K_X+B+\Mm_X)|_W$,  $L_{W^0}:=L_W|_{W^0}$, and $L:=L_W|_{W^n}$, where $W^n$ is the normalization of $W$. By induction hypothesis, $L$ is semi-ample$/U$.

Recall that $m>0$ is a sufficiently divisible integer such that $m(K_X+B+\Mm_X)$ is Cartier and $|m(K_X+B+\Mm_X)|_{X^0}|$ is base-point-free$/U^0$. Possibly replacing $m$ with a multiple, we may assume that 
\begin{itemize}
    \item $mL$ defines a contraction$/U$ $g^n: W^n\rightarrow Y^n$ such that there exists a very ample$/U$ divisor $H$ on $Y^n$ such that $(g^n)^*H=M$, and
    \item $mL_{W^0}$ defines a contraction$/U^0$ $g^0: W^0\rightarrow Z^0$, and there exists a very ample$/U^0$ divisor $H_{Z^0}$ on $Z^0$ such that $(g^0)^*H_{Z^0}=mL_{W^0}$.
\end{itemize}
In particular, one can checks that all conditions of Constructions \ref{cons: gluing part 1} and \ref{cons: glue part 2} hold. Therefore, in the following, we will adopt all notations as in Constructions \ref{cons: gluing part 1} and \ref{cons: glue part 2} (except that ``$\Delta$" will be replaced by ``$B$"). By Lemma \ref{lem: induced relation is stratified},
\begin{itemize}
    \item $(\sigma_1,\sigma_2):T^n\rightrightarrows Y^n$ and $(\sigma_{1H},\sigma_{2H}):T^n_H\rightrightarrows Y^n_H$ are stratified equivalence relations, and
    \item  $(Y^n,S_*(Y^n)),(T^n,S_*(T^n)),(Y^n_H,S_*(Y^n_H))$, $(T^n_H,S_*(T^n_H))$ satisfy (HN) and (HSN).
\end{itemize}
We let  $p_{Z^0}: Z^0_{H_{Z_0}}\rightarrow Z^0$ be the total spaces of the line bundle $H_{Z^0}$.

We let $Y^{n,0}=Y^n\times_UU^0$, $T^{n,0}=T^n\times_UU^0$, $Y^{n,0}_H=Y^n_H\times_UU^0$, and $T^{n,0}_H=T^n_H\times_UU^0$. Then the geometric quotients $Z^0=Y^{n,0}/T^{n,0}$ and $Z^0_{H_{Z^0}}=Y^{n,0}_H/T^{n,0}_{H}$ exist by \cite[Lemma 9.8]{Kol13}. In particular, the equivalence relations generated by $(\sigma_1,\sigma_2)|_{T^{n,0}}: T^{n,0}\rightrightarrows Y^{n,0}$ and $(\sigma_{1H},\sigma_{2H})|_{T^{n,0}_H}: T^{n,0}_H\rightrightarrows Y^{n,0}_H$ are finite. By \cite[Lemma 9.55]{Kol13}, the equivalence relations generated by $(\sigma_1,\sigma_2):T^n\rightrightarrows Y^n$ and $(\sigma_{1H},\sigma_{2H}):T^n_H\rightrightarrows Y^n_H$ are finite (cf. \cite[Proposition 3.12]{HX13}). By \cite[Theorem 9.21]{Kol13}, the geometric quotients $Y^n/T^n$ and $Y^n_H/T^n_H$ exist. 

We denote $Z:=Y^n/T^n$ and $Z_{H_Z}:=Y^n_H/T^n_H$. Then we have induced morphisms $p_Z: Z_{H_Z}\rightarrow Z$, $g: W\rightarrow Z$, and $\pi_Z: Y^n\rightarrow Z$, such that
\begin{itemize}
\item  $p_Z: Z_{H_Z}\rightarrow Z$ is a total space of a line bundle $H_Z$ on $Z$,
\item $Z^0=Z\times_UU^0$ and $Z^0_{H_{Z^0}}=Z_{H_Z}\times_UU^0$,
\item $g^0=g|_{W^0}$ and $g^*H_Z=mL_W$, and
\item $\pi_Z^*H_Z=H$.
\end{itemize}
Since $H$ is ample$/U$, $H_Z$ is ample$/U$. Thus $L_W$ is semi-ample$/U$. By Lemma \ref{lem: reduction to Nlc locus}, $K_X+B+\Mm_X$ is semi-ample$/U$, and we are done.
\end{proof}

The following theorem follows from Theorem \ref{thm: semi-ample over U0 implies semi-ample over U}.

\begin{thm}\label{thm: nef imply semi-ample}
Let $(X,B,\Mm)/U$ be an NQC gdlt g-pair and $A\geq 0$ an $\mathbb R$-Cartier $\mathbb R$-divisor on $X$. Assume that
\begin{enumerate}
    \item $K_X+B+\Mm_X$ is nef$/U$,
    \item $(X,B+A,\Mm)$ is glc, and
    \item $K_X+B+A+\Mm_X\sim_{\Rr,U}0$.
\end{enumerate}
Then $K_X+B+\Mm_X$ is semi-ample$/U$.
\end{thm}

\begin{proof}[Proof of Theorem \ref{thm: nef imply semi-ample}]
Possibly replacing $(X,B,\Mm)$ with a gdlt modification and replacing $A$ with the pullback of $A$, we may assume that $X$ is $\mathbb Q$-factorial. Since $-A$ is nef over $Z$, $\Supp A=f^{-1}(f(A))$. Since $(X,B+A,\Mm)$ is gdlt, $\Supp A$ does not contain any glc center of $(X,B,\Mm)$, hence $f(A)$ does not contain the image of any glc center of $(X,B,\Mm)$ in $U$. Let $U^0:=U\backslash f(A)$. Theorem \ref{thm: nef imply semi-ample} follows by applying Theorem \ref{thm: semi-ample over U0 implies semi-ample over U} to $(X,B,\Mm)/U$ and $U^0$ as $(K_X+B+\Mm_X)|_{X^0}\sim_{\mathbb R,U^0}0$, where $X^0:=X\times_UU^0$.
\end{proof}

\section{Du Bois property}\label{sec: DB singularity}
In this section we prove the g-pair versions of results in \cite[Chapter 6]{Kol13}, which will be used to prove Theorem \ref{thm: glc sings are Du Bois}. We adopt the notations as in \cite[Chapter 6]{Kol13} and will freely use them.

We first recall the following definition in \cite{Kov11} (cf. \cite[Definition 6.10]{Kol13}). 

\begin{defn}
A \emph{DB} pair $(X,\Sigma)$ consists of a reduced scheme $X$ of finite type and a closed reduced subscheme $\Sigma$ in $X$ such that the natural morphism
$$
\mathcal{I}_{\Sigma\subset X}\to \underline{\Omega}_{X,\Sigma}^0
$$
is a quasi-isomorphism. We will also say $(X,\Sigma)$ is \emph{DB} in this case.
\end{defn}

The definition of DB pairs is subtle but what really matters here is the following lemma:

\begin{lem}[{\cite[Proposition 6.15]{Kol13}}]\label{lem: property of DB pairs}
Let $(X,\Sigma)$ be a DB pair. Then $X$ has Du Bois singularities if and only if $\Sigma$ has Du Bois singularities.
\end{lem}

The following theorems are analogues of \cite[Theorems 6.31, 6.33]{Kol13} for g-pairs and the proofs are similar. For the reader's convenience, we provide full proofs here.

\begin{thm}\label{thm: (Z,W) is DB for glc crepant log structure}
Let $f:(X,\Delta,\Mm)\to Z$ be an NQC glc crepant log structure and $W\subset X$ the union of glc centers of $f:(X,\Delta,\Mm)\to Z$ except $Z$. Then $(Z,W)$ is a DB pair.
\end{thm}

\begin{proof}
By \cite[Theorem 2.28]{HL21a}, we may assume that $(X,\Delta,\Mm)$ is a $\Qq$-g-pair. Let $(Z,\Delta_Z,\NN)/U$ be a glc $\mathbb Q$-g-pair induced by the canonical bundle formula$/U$ of $f: (X,\Delta,\Mm)\rightarrow Z$ (the generalized canonical bundle formula). By Lemma \ref{lem: glc centers come from cbf}, the glc centers of $(Z,\Delta_Z,\NN)$ are exactly the glc centers of $f: (X,\Delta,\Mm)\rightarrow Z$. Thus we can assume that $f$ is the identity and $(X,\Delta,\Mm)=(Z,\Delta_Z,\NN)$.

Let $g:Y\to X$ be a log resolution such that $\Mm$ descends to $Y$ and $F:=g^{-1}(W)_{\red}$ is an snc divisor. Let
$$
K_Y+\Delta_Y+\Mm_Y:=g^*(K_X+\Delta+\Mm_X).
$$
and $D:=\Delta_Y^{=1}$. Since $\Mm_Y$ is nef$/X$ and big$/X$, there exists $0\le\Delta'_Y\sim_{\Qq,X}\Mm_Y$ such that $(Y,\Delta_Y-D+\Delta'_Y)$ is sub-klt. Let $\bar{\Delta}_Y:=(\Delta_Y-D+\Delta'_Y)^{\ge0}$ and $E:=(\Delta_Y-D+\Delta'_Y)^{\le 0}$, then $\lf\bar{\Delta}_Y\rf=0$ and $E$ is exceptional over $X$. Possibly replacing $Y$ with a higher resolution, we may assume that $D+E+\bar{\Delta}_Y$ is snc.

Since $E-D\ge-F$, we have natural maps:
$$
g_*\Oo_Y(-F)\to Rg_*\Oo_Y(-F)\to Rg_*\Oo_Y(E-D).
$$
Since $E-D\sim_{\Qq,X}K_Y+\bar{\Delta}_Y$, by \cite[Theorem 10.41]{Kol13},
$$
Rg_*\Oo_Y(E-D)\simeq_{qis}\sum_{i}R^ig_*\Oo_Y(E-D)[i].
$$
Thus we get a morphism 
$$
g_*\Oo_Y(-F)\to Rg_*\Oo_Y(-F)\to Rg_*\Oo_Y(E-D)\to g_*\Oo_Y(E-D).
$$
Note that 
$$
g_*\Oo_Y(E-D)=g_*\Oo_Y(E-D)\cap g_*\Oo_Y(E)=g_*\Oo_Y(E-D)\cap g_*\Oo_Y=g_*\Oo_Y(-D).
$$
Since $D$ is reduced and $g(D)=W$, we have $g_*\Oo_Y(-D)=\mathcal{I}_W$, the ideal sheaf of $W$ in $Z=X$. Moreover, $g_*\Oo_Y(-F)=\mathcal{I}_W$ since $F$ is also reduced. Therefore, we get an isomorphism $\mathcal{I}_W=g_*\Oo_Y(-F)\to g_*\Oo_Y(E-D)$, which implies that 
$$
\rho: \mathcal{I}_W\simeq g_*\mathcal{I}_F\to Rg_*\mathcal{I}_F
$$
has a left inverse. Since $Y$ is smooth and $F$ is an snc divisor, we see that $(Y,F)$ is a DB pair, thus by \cite[Theorem 3.3]{Kov12} (cf. \cite[Theorem 6.27]{Kol13}), $(Z,W)$ is also a DB pair. \end{proof}

\begin{thm}\label{thm: of glc origin implies DB}
Let $(X,S_*)$ be a stratified scheme of NQC glc origin (Definition \ref{defn: of glc origin}). Then $X$ is Du Bois.
\end{thm}

\begin{proof}
We use induction on the dimension.

Let $\pi: (X^n,S^n_*)\to (X,S_*)$ denote the normalization. Let $B(X)\subset X$ and
$B(X^n)\subset X^n$ denote the corresponding boundaries. By \cite[9.15.1]{Kol13}, we have a universal push-out diagram
\begin{center}
$\xymatrix{
B(X^n)\ar@{^(->}[r]\ar@{->}[d] & X^n\ar@{->}[d]^{\pi}\\
 B(X)\ar@{^(->}[r]& X\\
}$
\end{center}
Notice that $B(X)$ and $B(X^n)$ are of glc origin by Lemma \ref{lem: glc stratification is of glc origin}, hence Du Bois by induction.

Since $\pi$ is finite, it follows that $R\pi_*\mathcal{I}_{B(X^n)\subset X^n}=\pi_*\mathcal{I}_{B(X^n)\subseteq X^n}$. Furthermore, $\pi_*\mathcal{I}_{B(X^n)\subseteq X^n}=\mathcal{I}_{B(X)\subseteq X}$ by \cite[Theorem 9.30]{Kol13}. By \cite[Theorem 3.3]{Kov12} and Lemma \ref{lem: property of DB pairs}, we only need to show that $X^n$ is Du Bois. By assumption, for each irreducible component $X_i^n\subset X^n$ there is an NQC glc crepant log structure $f_i:(Y_i,\Delta_i,\Mm)\to Z_i$ and a finite surjection $Z_i\to X_i^n$. By \cite[Corollary 2.5]{Kov99}, we only need to show that $Z_i$ is Du Bois for each $i$. Let $B(Z_i)\subset Z_i$ be the boundary of the glc stratification of $Z_i$. Then  $B(Z_i)$ is of NQC glc origin by Lemma \ref{lem: glc stratification is of glc origin}, hence Du Bois by induction. By Theorem \ref{thm: (Z,W) is DB for glc crepant log structure}, $(Z_i,B(Z_i))$ is a DB pair, hence $Z_i$ is Du Bois and we are done.
\end{proof}

\begin{rem}
After finishing the first draft of the paper, the authors note the results \cite[Theorems 1,12]{KK22} proving the Du Bois property of varieties $V\subset X$ such that $\mld(V,X,\Delta)\leq\lcg(\dim X)$ for some lc pair $(X,\Delta)$, where $\lcg(\dim X)$ is the $1$-gap of lc thresholds. With the methods established in Sections \ref{sec: gluing} and \ref{sec: DB singularity}, we may also prove the NQC g-pair versions of \cite[Theorems 1,12]{KK22} by using the same arguments as in \cite{KK22}. In fact, as mentioned in \cite[Proof of Theorems 1 and 12]{KK22}, a quasi-log structure \cite{Fuj17} version of \cite[Theorems 1,12]{KK22} is expected and is used implicitly in \cite[Proof of Proposition 16]{KK22}, while any qlc pair is always an NQC glc g-pair (cf. \cite[Remark 1.9]{Fuj22}).
\end{rem}

\section{Proof of the main theorems}\label{sec: proof of the main theorems}

In this section we prove the main theorems, which are consequences of Theorems  \ref{thm: semi-ample over U0 implies semi-ample over U}, \ref{thm: nef imply semi-ample} and \ref{thm: of glc origin implies DB}.

\begin{proof}[Proof of Theorem \ref{thm: gmm over U0 implies gmm over U}]
By \cite[Theorem 3.14]{HL21a}, possibly replacing $(X,B,\Mm)$ with a gdlt model, we may assume that $(X,B,\Mm)$ is $\mathbb Q$-factorial gdlt. We run a $(K_X+B+\Mm_X)$-MMP$/U$ with scaling of an ample divisor
$$(X,B,\Mm):=(X_0,B_0,\Mm)\dashrightarrow (X_1,B_1,\Mm)\dashrightarrow\dots\dashrightarrow (X_i,B_i,\Mm)\dashrightarrow\dots.$$
By \cite[Lemma 2.7]{LX22}, possibly replacing $(X,B,\Mm)$ with $(X_n,B_n,\Mm)$ for some $n\gg 0$, we may assume that this MMP is an isomorphism over $U^0$ and $(X^0,B^0,\Mm^0)/U^0$ is a good minimal model of itself. Since every glc center of $(X,B,\Mm)$ intersects $X^0$ and $K_{X_i}+B_i+\Mm_{X_i}$ is semi-ample over $U^0$, $(X_i,B_i,\Mm)$ is log abundant$/U$ for any $i$. By \cite[Theorem 7.6]{LX22}, the MMP terminates with a log minimal model $(\bar X,\bar B,\Mm)/U$ of $(X,B,\Mm)/U$. By Theorem \ref{thm: semi-ample over U0 implies semi-ample over U}, $(\bar X,\bar B,\Mm)/U$ is a good minimal model of $(X,B,\Mm)/U$.
\end{proof}

\begin{proof}[Proof of Theorem \ref{thm: gmm exists for g-crepant log structure}]
 Since termination and semi-ampleness$/Z$ are both local on $Z$, we may assume that $Z$ is affine. By \cite[Theorem 1.3]{LX22} we get (1)(3). Possibly replacing $(X,B,\Mm)$ with $(Y,B_Y,\Mm)$ and replacing $A$ accordingly, we may assume that $K_X+B+\Mm_X$ is nef$/Z$, and we only need to show that $K_X+B+\Mm_X$ is semi-ample$/Z$.

Let $g: W\rightarrow X$ be a gdlt modification of $(X,B+A,\Mm)$, $0<\epsilon\ll 1$ a real number, $A_W:=g^*A$, and $K_W+B_W+A_W+\Mm_W:=g^*(K_X+B+A+\Mm_X)$, then $(W,\Delta_W:=B_W+(1-\epsilon)A_W,\Mm)$ is gdlt, $(W,\Delta_W+\epsilon A_W,\Mm)$ is glc, and $K_W+\Delta_W+\Mm_W+\epsilon A_W\sim_{\mathbb R,Z}0$. By Theorem \ref{thm: nef imply semi-ample}, 
$$\epsilon g^*(K_X+B+\Mm_X)\sim_{\mathbb R,Z}-\epsilon g^*A=-\epsilon A_W\sim_{\mathbb R,Z}K_W+\Delta_W+\Mm_W$$ 
is semi-ample$/Z$, hence $K_X+B+\Mm_X$ is semi-ample$/Z$, and we are done.
\end{proof}

\begin{proof}[Proof of Theorem \ref{thm: glc flip exists}]
Since $-(K_X+B+\Mm_X)$ is ample$/Z$, there exists an $\mathbb R$-divisor $0\leq A\sim_{\mathbb R,Z}-(K_X+B+\Mm_X)$ such that $(X,B+A,\Mm)$ is glc. By Theorem \ref{thm: gmm exists for g-crepant log structure}, there exists a good minimal model $(X',B',\Mm)/Z$ of $(X,B,\Mm)/Z$. We let $h: X'\rightarrow X^+$ be the birational morphism$/Z$ defined by $K_{X'}+B'+\Mm_{X'}$ and let $B^+:=h_*B'$.

We only need to show that the induce birational map $f^+: X^+\rightarrow Z$ is small. Let $p:W\to X$ and $q:W\to X'$ be a resolution of indeterminacy of $X\dashrightarrow X'$. Then $p^*(K_X+B+\Mm_X)=q^*(K_{X'}+B'+\Mm_{X'})+F$ where $F\geq 0$ is exceptional over $X'$. Let $D$ be a prime divisor on $X'$ that is exceptional over $X$, and $D_W$ its strict transform on $W$. Then $D_W$ is covered by a family of $p$-vertical curves $\Sigma _t$ such that $\Sigma_t\cdot p^*(K_X+B_X+\Mm_X)=0$. Since $F\cdot\Sigma_t\geq 0$, $\Sigma_t\cdot q^*(K_{X'}+B'+\Mm_{X'})\leq 0$. 
Let $\Sigma'_t=q_*\Sigma _t$, then $\Sigma ' _t\cdot (K_{X'}+B'+\Mm_{X'})\leq 0$ so that $\Sigma '_t$ are contracted by $X'\to X^+$ and hence $D$ is also contracted. Thus $X\dashrightarrow X^+$ does not extract any divisor, and $f^+$ is a $(K_X+B+\Mm_X)$-flip.
\end{proof}

\begin{proof}[Proof of Theorem \ref{thm: extracting divisor over glc structure}]
Let $g: Y\rightarrow X$ be a log resolution of $(X,\Supp B)$ such that $\Mm$ descends to $Y$ and $E$ is a divisor on $Y$. Let $a:=a(E,X,B,\Mm)\in[0,1)$ and $D:=\Supp\Exc(f)$. Let $B_Y:=g^{-1}_*B+D$, then $(Y,B_Y-aE,\Mm)$ is $\mathbb Q$-factorial gdlt. Thus $K_Y+B_Y-aE+\Mm_Y\sim_{\mathbb R,X}F\geq 0$ for some $\mathbb R$-divisor $F$ such that $E\not\subset\Supp F$. By \cite[Lemma 2.3]{LX22}, we may run a $(K_Y+B_Y-aE+\Mm_Y)$-MMP$/X$ with scaling of an ample divisor which terminates with a good minimal model $(W,B_W,\Mm)/X$ of $(Y,B_Y-aE,\Mm)/X$ and the induced birational map $Y\dashrightarrow W$ only contracts $F$. In particular, $E$ is still a divisor on $W$, and we let $E_W$ be the strict transform of $E$ on $E_W$. Then $\mult_{E_W}B_W=1-a>0$.

We may run a $(K_W+B_W-(1-a)E_W+\Mm_W)$-MMP$/X$ with scaling of an ample divisor. Since $(W,B_W,\Mm)/X$ is gdlt and $K_W+B_W+\Mm_W\sim_{\mathbb R,X}0$, by Theorem \ref{thm: gmm exists for g-crepant log structure}, this MMP terminates with a good minimal model $(Z',B_{Z'}-(1-a)E_{Z'},\Mm)/X$ of $(W,B_W-(1-a)E_W,\Mm)/X$, where $B_{Z'}$ and $E_{Z'}$ are the strict transforms of $B_W$ and $E_W$ on $Z'$ respectively. Thus  $-(1-a)E_{Z'}\sim_{\mathbb R,X}K_{Z'}+B_{Z'}-(1-a)E_{Z'}+\Mm_{Z'}$ is semi-ample$/X$, hence defines a birational morphism $Z'\rightarrow Z$ over $X$. We let $B_Z$ and $E_Z$ be the strict transforms of $B_{Z'}$ and $E_{Z'}$ on $Z$ respectively, and let $f: Z\rightarrow X$ be the induced morphism.

If $E_Z=0$, then $f$ is the identity map since $-E_Z$ is ample$/X$. Thus $B_Z=B$, so $a(E,Z,B_Z,\Mm)=a(E,X,B,\Mm)=a$. We have
\begin{align*}
    a&=a(E,Z,B_Z,\Mm)=a(E,Z',B_{Z'}-(1-a)E_{Z'},\Mm)\\
    &\geq a(E,W,B_W-(1-a)E_W,\Mm)>a(E,W,B_W,\Mm)=a,
\end{align*}
which is not impossible.

Therefore, $E_Z$ is a prime divisor on $Z$ and $-E_Z$ is ample over $X$, hence $\Supp E_Z$ contains all the exceptional locus on $Z$ and $f$ is an isomorphism away from $f(E_Z)$. In particular, $f$ only extracts $E$.
\end{proof}

\begin{proof}[Proof of Theorem \ref{thm: finite generation of R(X,A)}]
By \cite[Theorem 2.28]{HL21a}, we can assume $(X,B,\Mm)$ is a $\Qq$-g-pair. Let $g: Y\rightarrow X$ be a $\mathbb Q$-factorial gdlt modification of $(X,B,\Mm)$ and $K_Y+B_Y+\Mm_Y:=g^*(K_X+B+\Mm_X)$. Since $\Supp D$ does not contain any glc center of $(X,B,\Mm)$, then Cartier locus of $\Oo_X(-D)$ contains every generic point of the glc centers of $(X,B,\Mm)$. We may replace $D$ with $-A$ such that $A\ge 0$ and $\Supp A$ contains no glc center of $(X,B,\Mm)$. Let $0\leq C\sim -A$ be a divisor such that $C$ contains no glc centers of $(X,B,\Mm)$, then $A+C$ is Cartier and also contains no glc centers of $(X,B,\Mm)$. We may find an integral divisor $A_Y\le g^*(A+C)$ such that $A_Y\geq 0$ and $g(A_Y)=A$. 

Let $0<\epsilon\ll 1$ be a rational number and  $\Delta_Y:=B_Y+\epsilon g^*(A+C)-\epsilon A_Y$. Then $(Y,\Delta_Y+\epsilon A_Y,\Mm)$ is $\mathbb Q$-factorial gdlt and $K_Y+\Delta_Y+\epsilon A_Y+\Mm_Y\sim_{\Qq,X}0$. 
By Theorem \ref{thm: gmm exists for g-crepant log structure}, we may run a $(K_Y+\Delta_Y+\Mm_Y)$-MMP$/X$ which terminates with a good minimal model $(Z,\Delta_Z,\Mm)/X$ of $(Y,\Delta_Y,\Mm)/X$ with induced birational morphism $h: Z\rightarrow X$. Let $A_Z$ be the strict transform of $A_Y$ on $Z$, then $-A_Z\sim_{\Qq,X}K_Z+\Delta_Z+\Mm_Z$ is semi-ample$/X$, hence $R(Z/X,-A_Z)=R(X,-A)$ is a finite generated $\Oo_X$-algebra.
\end{proof}

\begin{proof}[Proof of Theorem \ref{thm: glc sings are Du Bois}]
Let $W$ be any union of the glc centers, then by Lemma \ref{lem: glc stratification is of glc origin} the induced stratified space $(W,S_*)$ is of NQC glc origin. Theorem \ref{thm: glc sings are Du Bois} follows from Theorem \ref{thm: of glc origin implies DB}.
\end{proof}

\end{document}